\documentclass[aos, seceqn, secthm,dvips]{arximspdf}
\usepackage{dcolumn}
\usepackage{graphicx}

\doi{10.1214/07-AOS566}
\volume{38}
\issue{3}
\pubyear{2010}
\firstpage{1885}
\lastpage{1912}

\makeatletter
\newcolumntype{d}[1]{D{.}{.}{#1}}
\newtheorem{prop}{Proposition}[section]
\newproclaim{remq}{Remark}[section]
\newtheorem{theorem}{Theorem}[section]
\newtheorem{cor}{Corollary}[section]
\newtheorem{lemma}{Lemma}[section]
\newcommand{\mano}{\matrix{\vdots\vspace*{-4pt}\cr\vdots\vspace*{-7pt}\cr\rule{1pt}{3pt}}}
\makeatother

\begin{document}
\begin{frontmatter}

\title{Efficient estimation for a subclass of shape invariant models}
\runtitle{Efficient estimation for shape invariant models}

\begin{aug}
\author[A]{\fnms{Myriam} \snm{Vimond}\corref{}\ead[label=e1]{myriam.vimond@ensai.fr}}
\runauthor{M. Vimond}
\affiliation{CREST-ENSAI, IRMAR}
\address[A]{Crest-Ensai\\
Campus de Ker Lann\\
Rue Blaise Pascal\\
BP 37203\\
35172 Bruz cedex\\
France\\
\printead{e1}\\
\phantom{E-mail:\ }} 
\end{aug}

\received{\smonth{5} \syear{2007}}
\revised{\smonth{8} \syear{2007}}

%
\begin{abstract}
In this paper, we observe a fixed number of unknown $2\pi$-periodic
functions differing
from each other by both phases and amplitude. This semiparametric model
appears in
literature under the name ``shape invariant model.'' While the common
shape is
unknown, we introduce an asymptotically efficient estimator of the
finite-dimensional
parameter (phases and amplitude) using the profile likelihood and the
Fourier basis.
Moreover, this estimation method leads to a consistent and
asymptotically linear
estimator for the common shape.
\end{abstract}

%
\end{frontmatter}

\section{Introduction}
In many studies, the response of interest is not a random variable but
a noisy function for each experimental unit, resulting in a sample of curves.
In such studies, it is often adequate to assume that the data
$Y_{i,j}$, the $i$th observation on the $j$th experimental unit,
satisfies the regression model
%
\begin{equation}\label{Model:general}
Y_{i,j}
=f^*_j(t_{i,j})+\sigma_j^*\varepsilon_{i,j},\qquad i=1, \ldots, n_j,
j=1, \ldots, J.
\end{equation}
Here, the unknown regression functions $f_j^*$ are $2\pi$-periodic and
may depend nonlinearly on the known regressors $t_{i,j}\in[0,2\pi]$.
The unknown error terms $\sigma_j^*\varepsilon_{i,j}$ are independent
zero mean random variables with variance ${\sigma_j^*}^2$.

The sample of individual regression curves will show a certain
homogeneity in structure, in the sense that curves coincide if they are
properly scaled and shifted.
In other words, the structure would be represented by the nonlinear
mathematical model
%
\begin{equation}\label{Model:shape:function}
f^*_j(t)=a_j^*f^*(t-\theta_j^*)+\upsilon^*_j\qquad\forall t \in
\mathbb{R}, \forall j=1, \ldots, J,
\end{equation}
where the shift $\theta^* = (\theta_j^*)_{j = 1, \ldots,  J}$,
the scale $a^* =  (a_j^*)_{j = 1, \ldots,  J}$ and the level
$\upsilon^* = (\upsilon_j^*)_{j = 1, \ldots,  J}$ are vectors of
$\mathbb R^J$ and the function $f^*$ is $2\pi$-periodic.
This semiparametric model was introduced by Lawton, Sylvestre and Maggio~\cite{Lawton} under the
name of shape invariant model.
We have both a finite-dimensional parameter
($\theta^*,a^*,\upsilon^*$) and an infinite-dimensional nuisance
parameter $f^*$ which is a member of some given large set of functions.
A general feature of semiparametric methods is to ``eliminate'' the
nonparametric component $f^*$, thus reducing the original
semiparametric problem to a suitably chosen parametric one.

Such models have been used to study child growth curves
(see \cite{Kneip}) or to improve a forecasting methodology
\cite{Loubes05} based on speed data of vehicles on a main trunk road
(see \cite{Gamboa04} for more details). Since the common shape is assumed to
be periodic, the model is particularly well adapted for the study of
circadian rhythms (see \cite{Wang}). Our model and our estimation
method are illustrated with the daily temperature of several cities.

The main goal of this paper is to present a method for the efficient
estimation of the parameter ($\theta^*,a^*,\upsilon^*$) without
knowing $f^*.$
The question of estimation of parameters for the shape invariant model
was studied by several authors.
First, Lawton, Sylvestre and Maggio \cite{Lawton} proposed an empirical procedure, \mbox{SEMOR}, based on
polynomial approximation of the common shape $f^*$ on a compact set.
The convergence and the consistency for SEMOR was proved by Kneip and Gasser~\cite{Kneip}.
H\"{a}rdle and Marron \cite{Hardle} built a $\sqrt n$~-consistent estimator and an
asymptotically normal estimator using a kernel estimator for the
function $f^*.$
Similar to Guardabasso, Rodbard and Munson \cite{Guardabasso}, Wang and Brown \cite{Wang}
and Luan and Li \cite{Genes}
used a smoothing spline
for the estimation of $f^*.$
The method of Gamboa, Loubes and Maza~\cite{Gamboa04} provides a $\sqrt n$-consistent
estimator and an asymptotically normal estimator for the shift
parameter~$\theta^*.$ This procedure is based on the discrete Fourier
transform of data. Our estimation method is related to the method of
Gamboa, Loubes and Maza \cite{Gamboa04}: The common shape $f^*$ is approximated by
trigonometric polynomials.

The efficiency of the estimators is to be understood as asymptotic
unbiasedness and minimum variance. To avoid the phenomena of
super-efficiency (e.g., Hodges estimators), the efficiency is
studied in a local asymptotic sense, under the local asymptotic
normality (LAN) structure.
The usual approach for determining the efficiency is to specify a least
favorable parametric submodel of the full semiparametric model (it is a
submodel for which the Fisher information is the smallest), locally in
a neighborhood of $f^*$, and to estimate ($\theta^*,a^*,\upsilon^*$)
in such a model (see \cite{VandVaart,VandVaart-StFlour}).
Here, we consider the parametric submodel where $f^*$ is a
trigonometric polynomial.
The method which is used is close to the procedure of
Gassiat and L{\'e}vy-Leduc \cite{Levy-Leduc} where the authors estimate efficiently the period of an
unknown periodic function.
The profile log-likelihood is used in order to ``eliminate'' the
nuisance parameter and to build an $M$-estimation criterion.
Moreover the efficiency of the $M$-estimator of
($\theta^*,a^*,\upsilon^*$) is proved by using the theory developed by
McNeney and Wellner \cite{Wellner}: The authors develop tools for nonindependent identically
distributed data that are similar in spirit to those for independent
identically distributed data.
Thus the notions of tangent space and of differentiability of the
parameter ($\theta^*,a^*,\upsilon^*$) are used in order to specify
the characteristics of an efficient estimator.
Under the assumptions listed in Theorem~\ref{thm:efficiency}, the
estimator of $(\theta^*,a^*,\upsilon^*)$ is asymptotically efficient.
This follows the conclusions of Murphy and Van der Vaart \cite{VandVaart-profile}:
Semiparemetric profile likelihoods, where the nuisance parameter has
been profiled out, behave like ordinary likelihoods in that they have a
quadratic expansion.

The profile log-likelihood induces the definition of an estimator for
the common shape.
Corollary~\ref{cor:shape:rate} establishes the consistency of this estimator.
The rate of the regression function estimator is the optimal rate in
nonparametric estimation \cite{VandVaart}, Chapter 24.
Using the theory developed by McNeney and Wellner \cite{Wellner}, we discuss its
efficiency: the estimator is asymptotically linear.
But the Fourier coefficients' estimators are efficient if and only if
the common shape $f^*$ is odd or even.
Even if this condition is satisfied, we can not deduce that the
estimator of $f^*$ is efficient because it is not regular.

This work is related to \cite{Empirical}, Chapter~3, where we propose another
criterion which allows us to estimate efficiently the parameter
($\theta^*,a^*,\upsilon^*$).
This criterion, which is similar by its definition to the criterion
proposed by Gamboa, Loubes and Maza \cite{Gamboa04} and~\cite{Empirical}, Chapter~2, allows us to build
a test procedure for the model.

The rest of the paper is organized as follows:
Section~\ref{sec:description} describes the model and the estimation method.
In Section~\ref{sec:trecia}, we discuss the efficiency of the estimator.
All technical lemmas and proofs are in Section~\ref{sec:ketvirta}.

\section{The estimation method}\label{sec:description}
\subsection*{The description of the model}
The data $(Y_{i,j})$ are the observations of $J$ curves at
the observation times $(t_{i,j}).$
We assume that each curve is observed at the same set of equidistant points
\[
t_{i}=t_{i,j}=\frac{i-1}{n} 2\pi\in[0,2\pi[,\qquad i=1, \ldots ,n.
\]
The choice of the observation times $t_i$ is related with the choice of
quadrature formula (see Remark~\ref{remq:quadrature:formula}).
The studied model is
%
\begin{equation}\label{Model}
\hspace*{6pt}Y_{i,j}=a_j^* f^*(t_{i}-\theta_{j}^*)+\upsilon_j^*
+\sigma^*\varepsilon_{i,j},\qquad j=1, \ldots ,J, i=1, \ldots ,n.
\end{equation}
The common shape $f^*$ is an unknown real $2\pi$-periodic continuous
function. We denote by $\mathcal{F}$ the set of $2\pi$-periodic
continuous functions. The noises ($\varepsilon_{i,j}$) are independent
standard Gaussian random variables.
For the sake of simplicity, we get a common variance
${\sigma^*}^2 = {\sigma_j^*}^2, j = 1, \ldots,  J.$
However, all our results are still valid for a general variance.

The model is semiparametric:
$\alpha^* = (\theta^*,a^*,\upsilon^*,\sigma^*)$ is the finite-dimensional parameter and $f^*$ is the
nuisance parameter.
Our aim is to estimate efficiently the internal shift
$\theta^* =(\theta_{j}^*)_{j = 1, \ldots,  J}$, the scale parameter
$a^* =(a_j^*)_{j = 1, \ldots,  J}$ and the external shift
$\upsilon^* = (\upsilon_{j}^*)_{j = 1, \ldots,  J}$ without knowing either the
shape $f^*$ or the noise level $\sigma^*$. We denote
$\mathcal{A}=[0,2\pi]^J\times\mathbb{R}^J\times[-\upsilon_{\max},\upsilon_{\max}]^J$
as the set where the parameter $(\theta^*,a^*,\upsilon^*)$ lies.

\subsection*{The identifiability constraints}
Before considering the estimation of parameters, we have to study the
uniqueness of their definition. Indeed, the shape invariant model has
some inherent unidentifiability: for a given parameter
$(\theta_0,a_0,\upsilon_0)\in\mathbb{R}^3$ and a shape function $f_0$ we can always
find another parameter $(\theta_1,a_1,\upsilon_1)\in\mathbb{R}^3$ and
another shape function $f_1$ such that
$a_0f_0(t-\theta_0)+\upsilon_0=a_1f_1(t-\theta_1)+\upsilon_1$ holds for all $t$.

Then we assume that the true parameters lie in the following spaces:
%
\begin{eqnarray}
\nonumber f^*\in\mathcal{F}_0
&=& \biggl\{f\in\mathcal{F}, c_0(f)=\int_0^{2\pi}f(t)\frac{dt}{2\pi}=0\biggr\}
\quad\mbox{and}\quad(\theta^*,a^*,\upsilon^*) \in \mathcal{A}_0,
\\
\eqntext{\mbox{where }\mathcal{A}_0 =
\Biggl\{(\theta,a,\upsilon)\in\mathcal{A}, \theta_1 = 0,
\displaystyle\sum_{j=1}^J  a_j^2 = J \mbox{ and } a_1>0
\Biggr\}.}
\end{eqnarray}
The constraint on the common shape allows us to uniquely define the
parameter $\upsilon^*\,[\upsilon_j^*=c_0(f_j^*)$, $j=1, \ldots, J$] and
to build asymptotically independent estimators (see Remark~\ref{rq:identifiability:constraint}).
The constant $\upsilon_{\mathrm{max}}$ is a user-defined (strictly positive)
parameter which reflects our prior knowledge on the level parameter.
The constraints $\theta_1=0$ and $a_1>0$ mean that the first unit
($j=1$) is taken as ``reference'' to estimate the shift parameter and
the scale parameter.
At last, the constraint $\sum_{j=1}^J  a_j^2 = J$ means that the
common shape is defined as the weighted sum of the regression functions
$f_j^*$ (\ref{Model:general}). This condition is well adapted to our
estimation criterion (see the next paragraph on the profile likelihood).

\subsection*{The profile log-likelihood}
Maximizing the likelihood function directly is not possible for
higher-dimensional parameters, and fails particularly for
semiparametric models.
Frequently, this problem is overcome by using a profile likelihood
rather than a full likelihood.
If $l_n(\alpha,f)$ is the full $\log$-likelihood, then the profile
likelihood for $\alpha\in\mathcal{A}_0$ is defined as
\[
pl_n(\alpha)=\sup_{f\in\mathcal{F}_0} l_n(\alpha,f).
\]
The maximum likelihood estimator for $\alpha$, the first component of
the pair $(\hat{\alpha}_n,\hat{f}_n)$ that maximizes
$l_n(\alpha,f)$, is the maximizer of the profile likelihood function
$\alpha\to pl_n(\alpha)$. Thus we maximize the likelihood in two steps. With the
assumptions on the model, we shall use the Gaussian $\log$-likelihood,
%
\begin{equation}\label{likelihood:eq}
l_n(\alpha,f)
=\frac{-1}{2\sigma^2}\sum_{i=1}^n \sum_{j=1}^J\bigl(Y_{i,j}-a_j f(t_i-\theta_j)-\nu_j \bigr)^2
-\frac{nJ}{2}\log{\sigma^2}.
\end{equation}

Generally, the problem of minimization on a large set is solved by the
consideration of a parametric subset.
Here, the semiparametric problem is reduced to a parametric one: $f$ is
approximated by its truncated Fourier series.
Thus the profile likelihood is approximated by minimizing the
likelihood $l_n$ on a subset of trigonometric polynomials.
More precisely, let $(m_n)_n$ be an increasing integer's sequence, and
let $\mathcal{F}_{0,n}$ be the subspace of $\mathcal{F}_0$ of trigonometric
polynomials whose degree is less than $m_n.$
In order to preserve the orthogonality of the discrete Fourier basis,
\begin{eqnarray*}
\forall|l| < \frac{n}{2}, \forall|p| < \frac{n}{2}\qquad\frac{1}{n}\sum_{r=1}^ne^{i(l-p) t_r}=
\cases{
1,&\quad if $l=p$,\cr
0,&\quad if $l\neq p$,
}
\end{eqnarray*}
we choose $m_n$ and $n$ such that
%
\begin{equation}\label{hyp:m_n}
2|m_n|<n,\qquad\lim_{n\to+\infty} m_n=+\infty\quad\mbox{and}\quad n\mbox{ is odd.}
\end{equation}
After some computations, the likelihood maximum is reached in the space
$\mathcal{F}_{0,n}$ by the trigonometric polynomial
%
\begin{equation}\label{likelihood:polynome}
\hat{f}_{\alpha}(t)= \sum_{1\leq|l|\leq m_n}\hat{c}_l(\alpha)
e^{il t}\qquad\forall t\in\mathbb{R},
\end{equation}
where for $l\in\mathbb{Z},$ $1\leq|l|\leq m_n,$
%
\begin{equation}\label{likelihood:coefficients}
\hspace*{28pt}\hat{c}_l(\alpha) =  \Biggl(n\sum_{j=1}^J a_j^2 \Biggr)^{-1}
\sum_{j=1}^J a_j  \sum_{i=1}^n (Y_{i,j}-\upsilon_j )
e^{-il(t_i-\theta_j)}\qquad\forall\alpha \in \mathcal{A}_0\times\mathbb{R}_+^*.
\end{equation}
Finally, using the orthogonality of the discrete Fourier basis, the
following equality holds:
\begin{eqnarray*}
&&\sum_{j=1}^{J}\sum_{i=1}^{n} \Biggl(Y_{i,j}-a_j\sum_{1\leq|l|
\leq m_n}\hat{c}_l(\alpha)e^{il(t_i-\theta_j)}-\nu_j \Biggr)^2
\\
&&\qquad=\sum_{j=1}^{J} \sum_{i=1}^{n} (Y_{i,j} - \nu_j)^2  -
\Biggl(n \sum_{j=1}^J a_j^2 \Biggr) \sum_{1\leq|l|<m_n}|\hat{c}_l(\alpha)|^2
\\
&&\qquad\quad{}+  n\sum_{1\leq|l|,|p|<m_n,l\neq p}
\hat{c}_l(\alpha)\overline{\hat{c}_p(\alpha)}\varphi_n
\biggl(\frac{l - p}{n}\biggr)
\sum_{j=1}^J  a_j^2e^{(p -
l)\theta_j},
\end{eqnarray*}
where $\varphi_n(t)=\sum_{s=1}^n e^{2i\pi st}/n.$
Let $M_n$ be the function of $\alpha=(\theta,a,\nu)$ defined as
\[
M_n(\alpha)= \frac{1}{nJ}   \sum_{j=1}^{J}\sum_{i=1}^{n}(Y_{i,j}
- \nu_j)^2-\sum_{1\leq|l|\leq m_n}  |\hat{c}_l(\alpha) |^2.
\]
With the identifiability constraints of the model, the profile
log-likelihood $pl_n$ is equal to
%
\begin{equation}\label{likelihood:criterion}
pl_n(\alpha)=- (nJ ) \frac{M_n(\alpha)}{2\sigma^2} -
\frac{nJ}{2}\log{\sigma^2}.
\end{equation}

\begin{remq}\label{remq:quadrature:formula}
The estimation method requires the estimation of the Fourier
coefficients of the common shape. A natural approach for estimating an
integral is to use a quadrature formula which is associated with the
observation times $t_i.$ In this paper, the observation times are
equidistant. Therefore the quadrature formula is the well-known
Newton--Cotes formula. Even if another choice of the observation times
is possible (see \cite{Empirical}, Chapter~2), this formula defines the discrete
Fourier coefficients $c_l^n(f)$ which are an accurate approximation of $c_l(f)$:
\[
c_l^n(f^*) =  \frac{1}{n} \sum_{s=1}^n   f^*(t_s) e^{-i l
t_s}\longrightarrow c_l(f^*) = \int_0^{2\pi}  f^*(t)e^{-ilt}\frac
{dt}{2\pi}.
\]
Moreover, the stochastic part of the coefficients (\ref{likelihood:coefficients}) are linear combinations of the complex
variables $w_{j,l},$
\[
w_{j,l}=\frac{1}{n}\sum_{r=1}^n e^{-ilt_r}\varepsilon_{i,j},\qquad j=1, \ldots, J, |l|\leq m_n.
\]
Due to Cochran's theorem, these variables are independent centered
complex Gaussian variables whose the variance is equal to $1/n.$ This
property is related to the convergence rate of the estimators (see
\cite{Empirical}, Chapter~2,
for more details, and \cite{Levy-Leduc} to compare).
\end{remq}

\subsection*{The estimation procedure}
Consequently, the maximum likelihood estimator of the finite-dimensional parameter is defined as
\[
\hat{\beta}_n=\arg\min_{\beta\in\mathcal{A}_0} M_n(\beta)\quad\mbox{or}\quad
\hat{\alpha}_n=(\hat{\beta}_n,\hat\sigma_n)=\arg\max_{\alpha\in\mathcal{A}_0\times\mathbb{R}_+^*} pl_n(\alpha).
\]
Then, the estimators of the common shape are the trigonometric
polynomials, which maximize the likelihood when $\alpha=\hat{\alpha}_n$:
\[
\hat{f}_n(t)=\hat{f}_{\hat{\alpha}_n}(t)= \sum_{1\leq|l|
\leq m_n}   \hat{c}_l(\hat{\alpha}_n) e^{il t}\qquad\forall t\in
\mathbb{R}.
\]
First, we study the consistency of the estimator of
($\theta^*,a^*,\upsilon^*$). The consistency of the common shape estimator is
studied in the next section.
\begin{theorem}[(Consistency)]\label{thm:consistance}
Assume that $2\pi$ is the minimal period of $f^*,$ and that
%
\begin{eqnarray}\label{thm:consist:assumption}
\sum_{|l| > m} |c_l(f^*)| =o \biggl(\frac{1}{\sqrt m} \biggr)\quad\mbox{and}\quad
\frac{m_n}{n}=o(1).
\end{eqnarray}
Then $\hat{\alpha}_n$ converges in probability to $\alpha^*.$
\end{theorem}

The assumption regarding the common shape means that the function $f^*$
is a $1/2$-holder function. The assumption on the number of Fourier
coefficients means that $m_n$ has to be small in relation to the number
of observation $n.$ Notice that Theorem \ref{thm:consistance} is still
valid if the noises $(\varepsilon_{i,j})$ are (centered) independent
identically distributed with finite variance.
\begin{pf*}{Proof of Theorem \ref{thm:consistance}}
The proof of this theorem follows the classical guidelines of the
convergence of $M$-estimators (see, e.g., Theorem 5.7 of Van der Vaart \cite{VandVaart}).
Indeed, to ensure consistency of $\hat{\beta}_n,$ it
suffices to show that:

\begin{longlist}
\item The uniform convergence of $M_n$ to a contrast function
$M+{\sigma^*}^2$ (Lem\-ma~\ref{lem:contraste}):
\[
\sup_{\beta\in\mathbb A} |M_n(\beta)-M(\beta)-{\sigma^*}^2 | =o_{P_\alpha^*}(1),
\]
where $M$ is defined as
\[
M(\beta)=\int_0^{2\pi}\frac{1}{J}\sum_{j=1}^J
\bigl(f_j^*(t)-\upsilon_j \bigr)^2\frac{dt}{2\pi}-\int_0^{2\pi}
\Biggl(\sum_{j=1}^J a_ja_j^*f^*(t-\theta_j^*+\theta_j) \Biggr)^2\frac
{dt}{2\pi}.
\]
\item $M(\cdot)$ has a unique minimum at $\beta^*$ (Lemma \ref{lem:minimum}).
\end{longlist}
\upqed
\end{pf*}

\subsection*{The daily temperatures of cities}
\begin{figure}[b]
\begin{tabular}{@{}c@{\quad}c@{}}
(a) & (b)
\\

\includegraphics{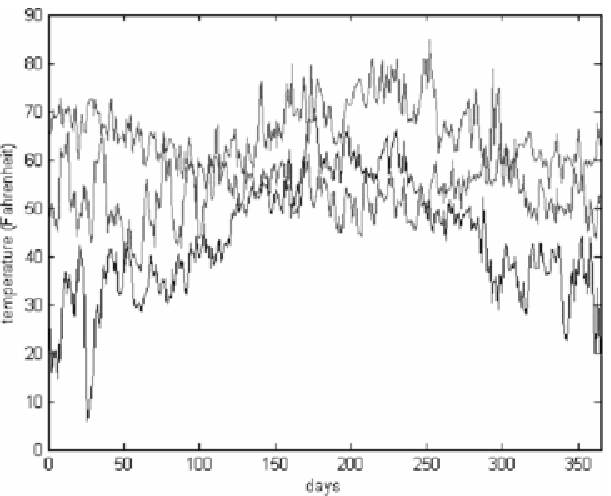}
 & \includegraphics{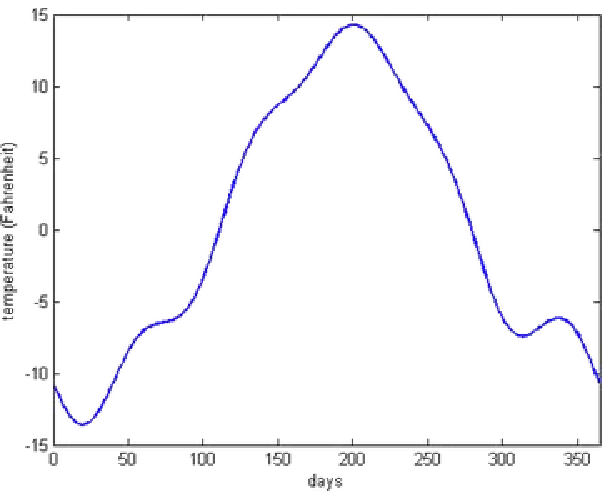}
\\
\end{tabular}
\caption{ \label{figtemperature}
(\textup{a}) Plots of the temperature curves associated with Juneau
(Alaska), Auckland (New Zealand) and Bilbao (Spain) in 2004.
(\textup{b}) Plot of the estimator of the common shape $\hat{f}_n$.
}
\end{figure}

The estimation method is applied to daily average temperatures (the
average daily temperatures are the average of 24 hourly temperature
readings). The data come from of the University of Dayton
(\url{http://www.engr.udayton.edu/weather/}). In order to
illustrate the method, we limit the study to three cities which have a
temperature range of an oceanic climate: Juneau (Alaska, city $j=1$),
Auckland (New Zealand, city $j=2$) and Bilbao (Spain, city $j=3$). An
oceanic climate is the climate typically found along the west coasts at
the middle latitudes of all the world's continents, and in southeastern
Australia. Similar climates are also found on coastal tropical
highlands and tropical coasts on the leeward sides of mountain ranges.
Figure \ref{figtemperature}(a) plots the sample of temperature curves.

If we assume that the data fit the model (\ref{Model}), the parameters
$\theta^*,$ $a^*$ and $\upsilon^*$ have the following meanings:
\begin{itemize}
\item $\upsilon_j^*$ is the annual temperature average of the $i$th city,
\item $a_j^*$ indicates whether the city is in the same hemisphere as
the first city ($a_j^*>0$) and measures the differences between the
winter and summer temperatures,
\item $\theta_j^*$ is the seasonal phase of the $i$th city,\vspace*{2pt}
\item $f^*$ describes the general behavior of the temperature evolution
of the oceanic climate.
\end{itemize}

The estimators of these parameters are given in Table \ref{tbl:estimators}.

\begin{table}
\caption{\label{tbl:estimators} Estimators of the parameters $\theta_2^*,$ $\theta_3^*,$
$a_1^*,$ $a_2^*,$ $a_3^*,$ $\upsilon_1^*,$ $\upsilon_2^*$ and
$\upsilon_3^*$}
\begin{tabular*}{\textwidth}{@{\extracolsep{\fill}}ld{2.4}d{2.4}d{2.5}@{}}
\hline
\textbf{City} $\bolds{j}$ & \multicolumn{1}{c}{$\bolds{j=1}$}
& \multicolumn{1}{c}{$\bolds{j=2}$} & \multicolumn{1}{c@{}}{$\bolds{j=3}$} \\
\hline
$\hat{\theta}_{j,n}$ (days)  & 0 & 12.5182 & 25.35381 \\
$\hat{a}_{j,n}$ & 1.2421 & -0.5833 & 1.0569 \\
$\hat{\upsilon}_{j,n}$ (Fahrenheit) & 43.9874 & 58.5312 & 60.1814\\
\hline
\end{tabular*}
\end{table}

Figure \ref{figtemperature}(b) plots the estimator of the common
shape. The number of the Fourier coefficients used to estimate the
common shape is $m_n=5.$ Further study will yield the most accurate
number $m_n,$ and leads to studying the estimation problem from the
point of view of the selection model.


\section{Efficient estimation}\label{sec:trecia}
\subsection{The LAN property}
Before studying the asymptotic efficiency of the estimators, we have to
establish the local asymptotic normality of the model. First, let us
introduce some notation. The model is semiparametric. The
finite-dimensional parameter $\alpha^*$ lies in
$\mathcal{A}_0\times\mathbb{R}_+^*.$ The nuisance parameter $f^*$ lies in $\mathcal{F}_0.$ For
$(\alpha,f)\in\mathcal{A}_0\times\mathbb{R}_+\times\mathcal{F}_0$ and
$t\in\mathbb R$,
we denote by $\mathbb{P}_{\alpha,f}(t)$ the Gaussian distribution in
$\mathbb{R}^J$ with variance $\sigma^2 I_J$ and mean
$(a_j f(t-\theta_j)+\nu_j)_{j=1,\ldots, J}$. Then the model of the observations is
\[
\mathcal{P}_n= \Biggl\{\mathbb{P}^{(n)}_{\alpha,f}
=\bigotimes_{i=1}^n \mathbb{P}_{(\alpha,f)}(t_i), (\alpha,f)\in
\mathcal{A}_0\times\mathbb{R}_+\times\mathcal{F}_0\Biggr\}.
\]

To avoid the phenomenon of super efficiency, we study the model on a
local neighborhood of $(\alpha^*,f^*).$ Let $(\alpha_n(h),f_n(h))$ be
close to $(\alpha^*,f^*)$ in the direction $h.$ The LAN property
requires that the $\log$-likelihood ratio for the two points
$(\alpha^*,f^*)$ and $(\alpha_n(h),f_n(h))$ converges in distribution to a
Gaussian variable which depends only on $h.$

Since the observations of our model are not identically distributed, we
shall follow the semiparametric analysis developed by McNeney and Wellner \cite{Wellner}.
The LAN property allows identification of the least favorable direction
$h$ that approaches the model, and thus allows us to know whether the
estimator is efficient. Let us denote the $\log$-likelihood ratio for
the two points $(\alpha^*,f^*)$ and $(\alpha,f)$
\[
\Lambda_n(\alpha,f)=\log\frac{d\mathbb{P}_{\alpha,f}^{(n)}}{
d\mathbb{P}_{\alpha^*,f^*}^{(n)}}.
\]

\begin{prop}[(LAN property)]\label{prop:LAN}
Assume that the function $f^*$ is not constant and is differentiable
with a continuous derivative denoted by $\partial f^*.$
Assume that the reals $a_j^*$, $j=1,\ldots, J,$ are nonnull. Considering
the vector space
$\mathcal{H}=\mathbb{R}^{J-1}\times\mathbb{R}^{J-1}
\times\mathbb{R}^{J}\times\mathbb{R}_+\times\mathcal{F}_0,$
the coordinates of a vector $h\in\mathcal{H}$ are denoted as follows:
\[
h=(h_{\theta,2},\ldots, h_{\theta,J},h_{a,2},\ldots,
h_{a,J},h_{\upsilon,2},\ldots,h_{\upsilon,J},h_{\sigma},h_f ).
\]
Then the space $\mathcal{H}$ is an inner-product space endowed with
the inner product $\langle\cdot,\cdot\rangle,$
\begin{eqnarray*}
\langle h,h\rangle
&=& J \frac{h_\sigma h_\sigma'}{{\sigma^*}^2}
+\frac{1}{{\sigma^*}^2} \Biggl\langle a_1^*h_f
- \sum_{j=2}^J h_{a,j}\frac{a_j^*}{a_1^*}f^*
\\
&&{}\hspace*{75pt}+ \upsilon_1,a_1^*{h'_f}- \sum_{j=2}^J {h'_{a,j}}\frac{a_j^*}{a_1^*}f^*
+ \upsilon_1'\Biggr\rangle_{\mathbb{L}^2}
\\
&&{}+\frac{1}{{\sigma^*}^2}\sum_{j=2}^J\langle a_1^*h_f + a_jf^* -
h_{\theta,j} a_j^* \partial  f^*
\\
&&\hspace*{55pt}{}+ h_{\upsilon,j},a_1^*{h'_f} +{h'_{a,j}}f^*
- {h'_{\theta,j}} a_j^* \partial  f^* +h_{\upsilon,j}'\rangle_{\mathbb{L}^2},
\end{eqnarray*}
where $\langle\cdot,\cdot\rangle_{\mathbb{L}^2}$ is the inner product in
$\mathbb{L}^2[0,2\pi].$
Moreover, the model (\ref{Model}) is LAN at $(\alpha^*,f^*)$ indexed
by the tangent space $\mathcal{H}.$
In other words, for each $h\in\mathcal{H},$ there exists a sequence
$(\alpha_n(h),f_n(h))$ such that
\[
\Lambda_n(\alpha_n(h),f_n(h))
=\Delta_{n}(h)-\tfrac{1}{2}\|h\|^2_{\mathcal{H}}+ o_{\mathbb{P}}(1).
\]
Here, the central sequence $\Delta_{n}(h)$ is linear with $h,$
\[
\Delta_{n}(h)=\frac{1}{\sqrt n}\sum_{i=1}^n\sum_{j=1}^J
\{(h_\sigma/\sigma^*)({\varepsilon^2_{i,j}}-1)+A_{i,j}^n(h) \varepsilon_{i,j}/\sigma^*  \} ,
\]
where for all $i=1,\ldots, n$,
\begin{eqnarray*}
A_{i,j}^n(h) =
\cases{
\displaystyle h_{a,1}^*f(t_i)-\sum_{k=2}^J h_{a,k} \frac{a_k^*}{a_1^*}f^*(t_i),&\quad if $j=1$,
\cr
a_j^*h_f(t_i-\theta_j^*)+h_{a,j}f^*(t_i-\theta_j^*)
\cr
\qquad{}-h_{\theta,j}a_j^*{\partial f^*(t_i-\theta_j^*)} + \nu_j,&\quad if $j = 2,\ldots, J$.
}
\end{eqnarray*}
\end{prop}

Notice that for the independent identically distributed semiparametric
models, the fact that the tangent space would not be complete does not
imply the existence of a least favorable direction.
In our model the tangent space $\mathcal H$ is a subset of the Hilbert space
\[
\overline{\mathcal H} = \mathbb R^{J -  1} \times
\mathbb R^{J -  1} \times \mathbb R^{J} \times \mathbb R \times
\{f\in \mathbb L^2[0,2\pi], c_0(f) = 0  \},
\]
endowed with the inner product $\langle\cdot,\cdot\rangle.$ Consequently, it is
easier to determine the least favorable direction using the Riesz
representation theorem.

\subsection{The efficiency}
The goal of this paper may be stated as the semiparametric efficient
estimation of the parameter
$\nu_n( \mathbb{P}_{\alpha^*,f^*}^{(n)}) = (\theta^*_2, \ldots, \theta^*_J,a^*_2
,\break\ldots, a_J^*,\upsilon^*_1 ,\ldots,\upsilon^*_J)$.
This parameter is differentiable relative to the tangent space
$\mathcal{H},$
\begin{eqnarray*}
&&\lim_{n\to\infty}  \sqrt n\bigl(  \nu_n \bigl( \mathbb{P}_{\alpha_n(h),f_n(h)}^{(n)}\bigr)
-  \nu_n \bigl( \mathbb{P}_{\alpha^*,f^*}^{(n)}  \bigr)  \bigr)\\
&&\qquad=  ( h_{\theta ,2} ,\ldots,  h_{\theta,J}, h_{a,2},\ldots ,h_{a,J},h_{\upsilon,2}, \ldots ,h_{\upsilon,J}) .
\end{eqnarray*}
Consequently, there exists a continuous linear map $\dot{\nu}$ from
$\mathcal H^{3J-2}$ on to $\mathbb{R}^{3J-2}.$
According to the Riesz representation theorem, there exist $3J-2$
vectors $(\dot{\nu}^{\theta}_j)_{2\leq j\leq J}$, $(\dot{\nu}^{a}_j)_{2\leq j\leq J}$
and $(\dot{\nu}^{\upsilon}_j)_{1\leq j\leq J}$
of $\overline{\mathcal H}$ such that
\[
\forall h \in \mathcal{H}\qquad \langle\dot{\nu}^{\theta}_j,h
\rangle=h_{\theta,j},\qquad \langle\dot{\nu}^{a}_j,h
\rangle=h_{a,j}\quad\mbox{and}\quad\langle\dot{\nu}^{\upsilon}_j,h
\rangle=h_{\upsilon,j}.
\]
These vectors are defined in Lemma \ref{lem:derivative:nu}.
Using the linearity with $h$ of $\Delta_{n}(h),$ the following
proposition, which is an application of Proposition 5.3 of McNeney and Wellner \cite{Wellner},
links the notion of asymptotic linearity of an estimator and
the efficiency.
\begin{prop}[(Asymptotic linearity and efficiency)]\label{prop:linearity}
Let $T_n$ be an asymptotically linear estimator of
$\nu_n(\mathbb{P}_{\alpha^*,f^*}^{(n)})$ with the central sequence
\[
\bigl(\Delta_{n}(\tilde{h}^{\theta}_2),\ldots,\Delta_{n}
(\tilde{h}^{\theta}_J),\Delta_{n}(\tilde{h}^{a}_2),\ldots,\Delta_{n}
(\tilde{h}^{a}_J),\ldots,\Delta_{n}(\tilde{h}^{\upsilon}_J)\bigr).
\]
$T_n$ is regular efficient if and only if for all $j$
$\tilde{h}^{\theta}_j=\dot{\nu}^{\theta}_j$, $\tilde{h}^{a}_j
=\dot{\nu}^{a}_j$ and $\tilde{h}^{\upsilon}_j=\dot{\nu}^{\upsilon}_j.$
\end{prop}

From Lemma~\ref{lem:derivative:nu}, if the assumptions of Proposition~\ref{prop:LAN} hold and if the estimator
$\hat{\beta}_n=(\hat\theta_n,\hat a_n,\hat\upsilon_n)$ is asymptotically linear, it is
efficient if and only if
\begin{eqnarray*}
\sqrt n (\hat\theta_n - \theta^* )
&=&\frac{\sigma^*}{\|\partial f^*\|_{\mathbb{L}^2}} \sum_{i=1}^n \biggl[
\frac{\mathbb{I}}{a_1^*}\mano\,{-}D^{-1}
\biggr]\partial F^*(t_i)\varepsilon_{i,\bolds{\cdot}}+o_{\mathbb{P}}(1),
\\
\sqrt n (\hat a_n - a^* )
&=&\frac{\sigma^*}{\| f^*\|_{\mathbb{L}^2}}
\sum_{i=1}^n \biggl\{\biggl[-\frac{a_1^*}{J}A \mano I_{J-1}-\frac{1}{J}A\,\, {}^t\!\!A\biggr]F^*(t_i)\biggr\}
\varepsilon_{i,\bolds{\cdot}} +o_{\mathbb{P}}(1),
\\
\sqrt n (\hat\upsilon_n - \upsilon^* )
&=&\sigma^*\sum_{i=1}^n\varepsilon_{i,\bolds{\cdot}}
+o_{\mathbb{P}}(1)\qquad\mbox{where }{}^t\varepsilon_{i,\bolds{\cdot}}
= {}^t(\varepsilon_{i1}, \ldots,\varepsilon_{i,J}),
\end{eqnarray*}
where $D$ is the diagonal matrix $\mathrm{diag}(a_2^* ,\ldots, a_J^*)$ and
$A={}^t(a_2^* ,\ldots, a_J^*)$ a vector in $\mathbb{R}^{J-1}.$ $F^*(t)$
and $\partial F^*(t)$ are, respectively, the diagonal matrix
$\mathrm{diag}(f^*\times(t-\theta_1^*) ,\ldots, f^*(t-\theta_J^*))$ and
$\mathrm{diag}(\partial f^*(t-\theta_1^*) ,\ldots,\partial f^*(t-\theta_J^*))$
for all $t\in\mathbb{R}$. We deduce the following theorem:
\begin{theorem}[(Efficiency)] \label{thm:efficiency}
Assume that the assumptions of Proposition \ref{prop:LAN} hold and that
%
\begin{eqnarray}
\sum_{l \in\mathbb{Z}} |l||c_l(f^*)| &<&\infty,\label{thm:efficiency:asump:f}\\
m_n^4/n&=&o (1). \label{thm:efficiency:asump:m}
\end{eqnarray}
Then $(\hat{\theta}_n,\hat{a}_n,\hat{\upsilon}_n)$ is
asymptotically efficient and
$\sqrt{n}(\hat{\theta}_n-{\theta}^*,\hat{a}_n-{a}^*,\hat{\upsilon}_n-{\upsilon}^*)$ converges in
distribution to a Gaussian vector $\mathcal{N}_{3J-2}(0,{\sigma^*}^2H^{-1}),$
where $H$ is the matrix defined as
\begin{eqnarray*}
H=
\pmatrix{
\|\partial f^*\|^2_{\mathbb{L}^2}  \biggl(D^2-\dfrac{1}{J}A^2\,\, {}^t\!\!A^2 \biggr) & 0 & 0
\cr
0 & \| f^*\|^2_{\mathbb{L}^2}  \biggl(I+\dfrac{1}{{a_1^*}^2}A\,\,  {}^t\!\!A \biggr) & 0
\cr
0 & 0 & I_J
}
\end{eqnarray*}
and its inverse matrix $H^{-1}$ is equal to
\begin{eqnarray*}
H^{-1}=
\pmatrix{
\dfrac{1}{\|\partial f^*\|^2_{\mathbb{L}^2}}\biggl(D^{-2}
+\dfrac{1}{{a_1^*}^2}\mathbb{I}_{J-1} {}^t\hspace*{-1pt}\mathbb{I}_{J-1} \biggr) & 0 & 0
\cr
0 & \dfrac{1}{\| f^*\|^2_{\mathbb{L}^2}} \biggl(I_{J-1} - \dfrac{1}{J}
A\,\, {}^t\!\!A \biggr) & 0
\cr
0 & 0 & I_J
}.
\end{eqnarray*}
\end{theorem}

\begin{pf}
Recall that the $M$-estimator is defined as the minimum of the
criterion function $M_n(\cdot).$ Hence, we get
\[
\nabla M_n(\hat{\beta}_n)=0,
\]
where $\nabla$ is the gradient operator.
Thanks to a second-order expansion, there exists $\bar{\beta}_n$ in a
neighborhood of $\beta^*$ such that
\[
\nabla^2 M_n(\bar{\beta}_n) \sqrt n (\hat{\beta}_n-\beta^*)=
-\sqrt n \nabla M_n(\beta^*),
\]
where $\nabla^2$ is the Hessian operator. Now, using two asymptotic
results from Proposition \ref{prop:gradient} and from Proposition~\ref{prop:hessian}, we obtain
\begin{eqnarray*}
\sqrt n (\hat\theta_n - \theta^* )
&=&\frac{\sigma^*}{\|\partial f^*\|^2_{\mathbb{L}^2}}\biggl(D^{-2}+\frac{1}{{a_1^*}^2}
\mathbb{I}_{J-1} {}^t\hspace*{-1pt}\mathbb{I}_{J-1}\biggr)G_n^\theta+o_{\mathbb{P}}(1),
\\
\sqrt n (\hat a_n - a^* )
&=&\frac{\sigma^*}{\| f^*\|^2_{\mathbb{L}^2}} \biggl(I_{J-1}
- \frac{1}{J} A\,\, {}^t\!\!A\biggr)G_n^a+o_{\mathbb{P}}(1),
\\
\sqrt n (\hat\upsilon_n - \upsilon^* )
&=&\sigma^*G_n^\upsilon+o_{\mathbb{P}}(1).
\end{eqnarray*}
\upqed
\end{pf}

\begin{figure}

\includegraphics{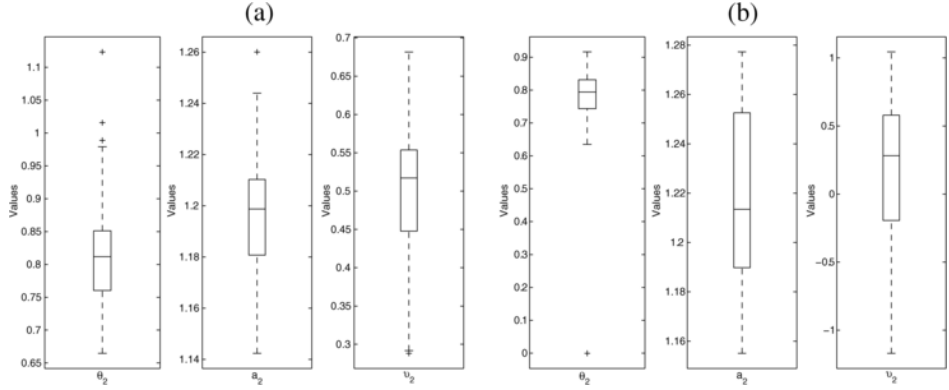}

  \caption{Boxplots of the estimators of $\theta_2^*$ $a_2^*$ and
$\upsilon_2^*$ associated with the space parameter $\mathcal{A}_0$
(\textup{a}) and
$\mathcal{A}_1$ (\textup{b}). The data are generated with $f^*(t)=20*t/(2\pi)(1-t/(2\pi))$,
$\theta^*=(0 \quad0.8)$, $a^*=( 0.75 \quad1.1990),$
$ \nu^*=(7.5/3 \quad0.5)$ and $n=201.$ The boxplots are computed from
$100$ sets of data.}\label{fig:boxplot}
\end{figure}

\begin{remq}\label{rq:identifiability:constraint}
The choice of the identifiability constraints is important for the
relevancy of the estimation. For example, if we no longer assume that
$c_0(f)$ is null, we may consider the following parameter space:
\[
\mathcal{A}_1
=\Biggl\{(\theta,a,\nu)\in\mathcal{A},
\mbox{ such that }\theta_1 = 0, \sum_{j=1}^J a_j^2
= J\mbox{ and }a_1 >0 \Biggr\}\quad\mbox{and}\quad f\in\mathcal{F}.
\]
Consequently we have to estimate $3J - 3$ parameters: $\theta_2^*,
\ldots ,\theta_J^*$, $a_2^*, \ldots ,a_J^*$, and $\upsilon_2^*,
\ldots ,\upsilon_J^*.$
This choice modifies the estimation criterion and the tangent space, too.
Nevertheless, if the assumptions of Theorem \ref{thm:efficiency} hold,
the estimator is asymptotically efficient. But its covariance matrix is
not block diagonal any more:\looseness=1
\begin{eqnarray*}
\Gamma = {\sigma^*}^2
{\fontsize{8}{10}\selectfont
\pmatrix{
\dfrac{1}{\|\partial f^*\|^2_{\mathbb{L}^2}}
\biggl(D^{-2} +\dfrac{1}{{a_1^*}^2}\mathbb{I}_{J-1} {}^t \mathbb{I}_{J - 1}\biggr) & 0 & 0
\cr
0 & \dfrac{1}{\| f^*\|^2_{\mathbb{L}^2} - c_0(f^*)^2} B
& \dfrac{-c_0(f^*)}{\| f^*\|^2_{\mathbb{L}^2} - c_0(f^*)^2}I_{J-1}
\cr
0 & \dfrac{-c_0(f^*)}{\| f^*\|^2_{\mathbb{L}^2} - c_0(f^*)^2}I_{J-1}
& \dfrac{\|f^*\|^2_{\mathbb{L}^2} }{\| f^*\|^2_{\mathbb{L}^2} - c_0(f^*)^2} B^{-1}
}}
,
\end{eqnarray*}
where $B=I_{J-1} - \frac{1}{J} A\,\, {}^t\!\!A$ with
$B^{-1}=I_{J-1} + \frac{1}{{a_1^*}^2} A\,\, {}^t\!\!A.$
In other words, $\hat a_n$ and $\hat\upsilon_n$ are not
asymptotically independent: modifying the identifiability constraint
$c_0(f^*)=0$ damages the quality of the estimation.

To illustrate this phenomenon, we present the boxplots of the
estimators which are relatively associated with the parameter space
$\mathcal{A}_0$ [Figure~\ref{fig:boxplot}(a)] and $\mathcal{A}_1$
[Figure~\ref{fig:boxplot}(b)]. Let $(\alpha^*,f^*)$ be a parameter of the model.
With the constraints associated with the parameter space $\mathcal{A}_0,$ we
have to estimate $\theta_2^*$, $a_2^*$, $\upsilon_1^*$ and $\upsilon_2^*$ for the following model ($J=2$):
\begin{eqnarray*}
\cases{
Y_{i,1}=a_1^*f^*(t_i)+\upsilon_1^*+\varepsilon_{i,1},&\quad$i=1,\ldots,n$,
\cr
Y_{i,2}=a_2^*f^*(t_i-\theta_2^*)+\upsilon_2^*+\varepsilon_{i,2},&\quad$i=1,\ldots,n$.
}
\end{eqnarray*}
With the constraints associated with the parameter space $\mathcal{A}_1,$ we
have to estimate $\theta_2^*$, $a_2^*$ and $\upsilon_2.$ The data may
be rewritten as
\begin{eqnarray*}
\cases{
Y_{i,1}=a_1^*g^*(t_i)+\varepsilon_{i,1},&\quad$i=1,\ldots,n$,
\cr
Y_{i,2}=a_2^*g^*(t_i-\theta_2^*)+\upsilon_2+\varepsilon_{i,2},&\quad$i=1,\ldots,n$,
}
\end{eqnarray*}
where $g^*=f^*+\upsilon_1^*$ and $\upsilon_2=\upsilon_2^*-a_2^*\upsilon_1^*.$
After generating several sets of data from a parameter
$(\alpha^*,f^*)$ which we have chosen, we have computed the estimators of
$\theta^*_2$, $a_2^*$ and $\upsilon^*_2$ for every set of data. Figure~\ref{fig:boxplot}
presents the boxplots of the estimators of
$\theta^*_2$, $a_2^*$ and $\upsilon^*_2$ for these two models.
\end{remq}

As a consequence of the previous theorem, the Gaussian vector $G_n$
converges in distribution to a centered Gaussian vector $\mathcal
N_{3J-2} (0,H ),$ and the equation holds:\looseness=1
\[
\sqrt n  (\hat\beta_n-\beta^* )= (H/{\sigma
^*}^2 )^{-1} \sigma^*G_n+o_{\mathbb{P}}(1).
\]
Comparing this formula with the results of the independent identically
distributed semiparametric model (see \cite{VandVaart}), we identify
the efficient information matrix as $H/{\sigma^*}^2$ and the efficient
score as $\sigma^*G_n.$

Indeed, let $X_1,\ldots, X_n$ be a random sample from a distribution
$\mathbb P$ that is known to belong to a set of probabilities $\{
\mathbb P_{\theta,\eta}, \theta\in\Theta\subseteq\mathbb R^d,\eta
\in\mathcal G\}$.
Then an estimator sequence $T_n$ is asymptotically efficient for
estimating $\theta$ if
\[
\sqrt n(T_n-\theta)= (\tilde I_{\theta,\eta} )^{-1}
\Biggl(\frac{1}{\sqrt n}\sum_{i=1}^n\tilde l_{\theta,\eta}(X_i)\Biggr)+o_{\mathbb{P}}(1),
\]
where $\tilde l_{\theta,\eta}$ is the efficient score function, and
$\tilde I_{\theta,\eta}$ is the efficient information matrix.

Moreover, our result follows Murphy and Van der Vaart \cite{VandVaart-profile}.
The authors demonstrate that if the entropy of the nuisance parameters
is not too large and the least favorable direction exists, the profile
likelihood behaves very much like the ordinary likelihood and the
profile likelihood correctly selects a least favorable direction for
the independent identically distributed semiparametric model.
This holds if the profile $\log$-likelihood $pl_n$ verifies the
following equation:
\begin{eqnarray*}
&&pl_n( \hat\theta_n )-pl_n(\theta)
\\
&&\qquad = \sum_{i=1}^n \tilde l_{\theta,\eta}(X_i)(\hat\theta_n -\theta)-\tfrac{1}{2}n\,{}^t\!(\hat\theta_n -\theta)\tilde I_{\theta,\eta}
(\hat\theta_n-\theta)+o_{\mathbb{P}}\bigl( \sqrt n\|\hat\theta_n -\theta\|+1 \bigr)^2
 ,
\end{eqnarray*}
where $\hat\theta_n$ maximizes $pl_n.$
Then, if $\tilde I_{\theta,\eta}$ is invertible, and $\hat\theta_n$
is consistent, $\hat\theta_n$ is asymptotically efficient.

For our model, a similar asymptotic expansion holds.
Indeed, by a Taylor expansion, there exists $\tilde\alpha_n$ such that
\begin{eqnarray*}
&&pl_n( \hat\alpha_n ) - pl_n( \alpha^* )
\\
&&\qquad= n^{1 / 2}G_n( \hat\beta_n - \beta^* ) - \frac{n}{2}\,\,{}^t\!( \hat\beta_n - \beta^* )\frac{H}{{\sigma^*}^2}
( \hat\beta_n - \beta^* )
\\
&&\qquad\quad+ o_{\mathbb{P}}
( n^{1 / 2}
\| \hat\beta_n - \beta^* \| + 1 )^2  .
\end{eqnarray*}
%

\subsection{Asymptotic linearity of the common shape estimator}\label{subsec:estimator:shape}
In this subsection, we study the consistency and the characteristics of
the estimator of the
common shape which is defined in Section \ref{sec:description}.
We show that the convergence rate of this estimator is the optimal rate
for the nonparametric estimation.

\begin{cor}\label{cor:shape:rate}
Assume that $f^*$ is $k$ times continuously differentiable with
$\int_0^{2\pi}  |f^{(k)} (t)|^2\,dt < \infty$ and $k\geq1.$
Furthermore, suppose that the assumptions of Theorem \ref{thm:efficiency} hold;
then there exists a constant $C$ such that for a
large $m_n$
\begin{eqnarray*}
\sup_{t\in\mathbb{R}} |\hat{f}_n(t)-f^*(t) |
&=& O_{\mathbb{P}} \biggl(\frac{1}{m_n^{k-1/2}}+\frac{m_n}{\sqrt n} \biggr),
\\
\int_0^{2\pi}\mathbb E  \bigl(\hat{f}_n(t)-f^*(t) \bigr)^2
&\leq&
C \biggl(\frac{1}{m_n^{2k}}+\frac{m_n}{n} \biggr).
\end{eqnarray*}
Consequently, for $m_n\sim n^{1/(2k+1)}$, we have $\mathit{MISE}_{f^*}(\hat
{f}_n)=O (n^{-2k/(2k+1)} ).$
\end{cor}

Let $\mathbb B$ represent the Banach space defined as the closure of
$\mathcal F$ for the $\mathbb{L}^2$-norm
\[
\mathbb B= \{f\in\mathbb{L}^2[0,2\pi] \mbox{ such that }
c_0(f)=0 \}.
\]
Here, the studied sequence of parameter $\nu_n$ is not $(\theta
^*,a^*,\upsilon^*)$ any more, but it is the truncated Fourier series
of $f^*$:
\[
\nu_n\bigl(\mathbb{P}_{\alpha^*,f^*}^{(n)}\bigr)
=\sum_{|l|\leq m_n}c_l(f^*)e^{il(\cdot)}.
\]
The parameter sequence $\nu_n$ is differentiable:
\[
\lim_{n\to\infty}\sqrt n \bigl( \nu_n
\bigl(\mathbb{P}_{\alpha_n(h),f_n(h)}^{(n)} \bigr) - \nu_n
\bigl(\mathbb{P}_{\alpha^*,f^*}^{(n)} \bigr) \bigr)=h_f.
\]
Thus, there exists a continuous linear map $\dot\nu$ from $\mathcal
H$ on to $\mathbb B.$ To have a representation of the derivative $\dot
\nu,$ we consider the dual space $\mathbb{B}^*$ of $\mathbb{B}.$
In other  words, for $b^*\in\mathbb{B}^*,$ $b^*\dot\nu$ is
represented by $\dot{\nu}^{b^*}\in\overline{\mathcal H}$:
\[
\forall h \in  \mathcal{H}\qquad b^*\dot\nu(h)= \langle\dot{\nu
}^{b^*},h \rangle=b^*h_f.
\]
Furthermore, the dual space $\mathbb{B}^*$ is generated by the
following linear real functions:\looseness=1
\begin{eqnarray*}
b^*_{1l}\dvtx f \in \mathcal{F}_0
&\to&\int_0^{2\pi}
f(t)\cos(lt)\frac{dt}{2\pi}\quad\mbox{and}
\\
\quad b^*_{2l}\dvtx f \in
\mathcal{F}_0
&\to& \int_0^{2\pi}  f(t)\sin(lt)\frac{dt}{2\pi},\qquad l\in\mathbb{Z}^*.
\end{eqnarray*}
Thus it suffices to know $\dot{\nu}^{b^*_{1l}}$ and $\dot{\nu}^{b^*_{2l}}$
for all $l\in\mathbb{Z}^*$ in order to determine all $\{\dot{\nu}_{b^*},b^*\in\mathbb{B}^*\}.$
After straightforward computations, these vectors are
\[
\dot{\nu}_{b^*_{1l}}= (0,\cos(l\cdot)/J )\quad\mbox{and}\quad\dot{\nu}_{b^*_{2l}}= (0,\sin(l\cdot)/J ).
\]

The estimator of the common shape is asymptotically linear. This means
that for all $b^*\in\mathbb{B}^*$ there exists $h^{b^*}\in\overline
{\mathcal{H}}$ such that
%
\begin{equation}\label{eq:asymptotic:linear}
\sqrt n b^* \bigl(T_n-\nu_n(\mathbb{P}_{\alpha^*,f^*}^n)\bigr)
=\Delta_n(h^{b^*})+o_{\mathbb P}(1).
\end{equation}
Since $ \{b^*_{1l},b^*_{2l}, l\in\mathbb{Z}^* \}$
generates the dual space of $\mathbb{B},$
Lemma \ref{lem:asymptotic:linear:f} ensures the asymptotic linearity
of $\hat{f}_n.$

Now, we discuss the regularity and the efficiency of this estimator.
We deduce from Proposition 5.4 of McNeney and Wellner \cite{Wellner} that:
\begin{cor}
$b^*\hat{f}_n$ is a regular efficient estimator of $b^*f^*$ for all
$b^*\in\mathbb{B}^*$ if and only if the function $f^*$ is odd or even.
In particular, in this case, the estimator of the Fourier coefficients
of $f^*$ is efficient.
\end{cor}

Consequently, $\hat{f}_n$ is eventually regular and efficient if the
common shape $f^*$ is odd or even. But the fluctuations
$
\sqrt n (T_n-\nu_n (\mathbb{P}_{\alpha_n(h),f_n(h)}^{(n)} )  )
$
do not converge weakly under
$
\mathbb{P}_{\alpha_n(h),f_n(h)}^{(n)}
$
to a tight limit in $\mathbb B$ for each
$ \{\alpha_n(h),f_n(h) \}$ [e.g., take $h=(0,0)$]. Thus, even if $f^*$ is
odd or even, $\hat{f}_n$ is not efficient.

\begin{remq}
The model where the function $f^*$ is assumed to be odd or even has
been studied by Dalalyan, Golubev and Tsybakov \cite{Dalalyan06}. In this model, the identifiability
constraint ``$\theta_1=0$'' is not necessary: The shift parameters are
defined from the symmetric point $0.$ Thus the estimator of
$\theta_1^*,\ldots,\theta_J^*$ would be asymptotically independent. Moreover
the estimation method would be adaptative.\looseness=1
\end{remq}

\section{The proofs}\label{sec:ketvirta}
\subsection{\texorpdfstring{Proof of Theorem
\protect\ref{thm:consistance}}{Proof of Theorem 2.1}}
\begin{remq}\label{remqnotations}
Let us introduce some notation. First the deterministic part of
$\hat{c}_l$ (\ref{likelihood:coefficients}) is equal to
%
\begin{eqnarray}
&&\frac{1}{nJ}\sum_{j=1}^J\sum_{i=1}^n a_ja_j^*f^*(t_i - \theta_j^*)
e^{il(t_i -\theta_j)}
\nonumber\\
&&\qquad= \sum_{p\in\mathbb{Z}}c_p(f^*)\varphi_n\biggl(\frac{l -p}{n} \biggr)
\phi(l\theta - p\theta^* ,a)\nonumber
\\
&&\qquad= c_l(f^*)\phi(l\theta - l\theta^*,a)+g_n^l(\beta)\label{egaliteg_n}
\\
&&\qquad\quad\eqntext{\displaystyle\mbox{where }g_n^l(\beta)
= \mathop{\sum_{|p|\geq m}}_{p-l \in n\mathbb{Z}} c_p(f^*)
\phi(l\theta - p\theta^* ,a)\quad\mbox{and}}
\\\nonumber
&&\qquad\quad\hspace*{48pt}\phi(\theta,a) = \sum_{j=1}^J   a_ja_j^*e^{i\theta_j}/J.
\end{eqnarray}
Since assumption (\ref{hyp:m_n}) holds, the term $g_n^l$ is bounded by
%
\begin{equation}\label{proof:borne:gnl}
 |g_n^l(\beta) |\leq\sum_{2|p|\geq n} |c_p(f^*)|.
\end{equation}
For $j=1,\ldots, J$ and $|l|\leq m_n,$ let us denote the variable
$\xi_{j,l}$ as $w_{j,l}=\xi_{j,l}/\sqrt n.$ Then the variables $\xi_{j,l}$
are independent standard complex Gaussian variables from Remark
\ref{remq:quadrature:formula}.
Thus the stochastic part of $\hat{c}_l$ is equal to
%
\begin{eqnarray}\label{proof:borne:xil}
\hspace*{10pt}\frac{\sigma^*}{\sqrt n}\xi_l(\beta)
&=& \frac{\sigma^*}{J\sqrt n}\sum_{j=1}^J a_je^{il\theta_j}\xi_{j,l}
\qquad\mbox{with }|\xi_l(\beta) |\leq \frac{\sigma^*}{J\sqrt n}\sum_{j=1}^J |\xi_{j,l}|.
\end{eqnarray}
\end{remq}

\begin{lemma}[(The uniform convergence in probability)]\label{lem:contraste}
Under the assumptions of Theorem \ref{thm:consistance}, we have
\[
\sup_{\beta\in\mathcal{A}_0} |M_n(\beta)-M(\beta)-{\sigma^*}^2 | =o_{P_\beta^*}(1),
\]
where $ M(\beta)=M^1(\beta)+M^2(\beta),$
\begin{eqnarray*}
M^1(\beta)=\sum_{l\in\mathbb{Z}^*} |c_l(f)|^2
\bigl(1-|\phi(l\theta-l\theta^*,a)|^2 \bigr)\quad\mbox{and}\quad M^2(\beta)
= \frac{1}{J}\sum_{j=1}^J (\upsilon_j^*-\upsilon j  )^2.
\end{eqnarray*}
\end{lemma}
\begin{pf}
The contrast process may rewritten as the sum of three terms:
\[
M_n(\beta)=D_n(\beta)+\sigma^*L_n(\beta)+{\sigma^*}^2 Q_n(\beta).
\]
The term $D_n(\beta)=D_n^1(\beta)-D_n^2(\beta)$ is the deterministic
part where
\begin{eqnarray*}
D_n^1(\beta)
&=& \frac{1}{Jn} \sum_{j=1}^J
\Biggl\{\sum_{i=1}^n a_j^*f^*(t_i - \theta_j^*) + \upsilon_j^* - \upsilon_j \Biggr\},
\\
D_n^2(\beta)
&=&  \sum_{1\leq|l|\leq m_n} \Biggl|\sum_{p\in
\mathbb{Z}}c_p(f^*)\varphi_n   \biggl(\frac{l - p}{n} \biggr)\phi
(l\theta - p\theta^* ,a)\Biggr|^2.
\end{eqnarray*}
The term $L_n(\beta)=L^1_n(\beta)-L^2_n(\beta)$ is the linear part
with noise, where
\begin{eqnarray*}
L^1_n(\beta)&=&\frac{2}{nJ}\sum_{j=1}^J\sum_{i=1}^n
\bigl(a_j^*f^*(t_i-\theta_j^*)+\upsilon_j^*-\upsilon_j \bigr)\sigma^*\varepsilon_{i,j},
\\
L^2_n(\beta)&=&\frac{2}{\sqrt n}\sum_{1\leq|l|\leq m_n} \Re \Biggl\{
\sum_{p\in\mathbb{Z}}c_p(f^*)\varphi_n \biggl(\frac{l-p}{n}\biggr)
\phi(l\theta-p\theta^*,a) \overline{\xi_l(\beta)} \Biggr\}.
\end{eqnarray*}
The term $Q_n(\beta)=Q^1_n(\beta)-Q^2_n(\beta)$ is the quadratic
part with noise:
\begin{eqnarray*}
Q^1_n(\beta)=\frac{1}{nJ}\sum_{j=1}^J\sum_{i=1}^n\varepsilon_{i,j}^2\quad\mbox{and}\quad Q^2_n(\beta)
=\frac{1}{n}\sum_{1\leq|l[<m_n}  |\xi_l(\beta) |^2.
\end{eqnarray*}

From the weak law of large numbers, $Q_n^1$ does not depend on $\beta$
and converges in probability to $1$. Furthermore, $Q_n^2$ is bounded by
\[
0\leq Q_n^2(\beta) \leq Q_n^B\qquad\mbox{where }nJ Q_n^B =
\sum_{|l|<m_n}  \sum_{j=1}^J    |\xi_{j,l} |^2.
\]
Then assumption (\ref{thm:consist:assumption}) induces that
$\sup_{\beta\in\mathcal{A}_0} |Q_n(\beta)-1 |$ converges to
$0$ in
probability.

Using the fact that $f^*$ is continuous and that
$|\upsilon_j|\leq\upsilon_{\max}$, there exists a constant $c>0$ such that for all
$\beta\in\mathcal{A}_0$ we have
\[
|L_n^1(\beta) |\leq c L_n^{1B}\qquad\mbox{where }L_n^{1B}
= \frac{1}{nJ} \Biggl|\sum_{j=1}^J\sum_{i=1}^n\varepsilon_{i,j} \Biggr|.
\]
Then we deduce that $L_n^1$ converges uniformly in probability to~$0$.
Concerning the term $L_n^2,$ it may be written as the sum of two
variables $L_n^{21}$ and $L_n^{22}$:
\begin{eqnarray*}
\sqrt n L_n^{21}(\beta)
&=&
2\Re \Biggl\{\sum_{1\leq|l|\leq m_n}
c_l(f^*)\phi(l\theta-l\theta^*,a)\overline{\xi_l(\beta)}\Biggr\},
\\
\sqrt n L_n^{22}(\beta)
&=&
2\Re \Biggl\{\sum_{1\leq|l|\leq m_n}
 g_n^l(\beta)\overline{\xi_l(\beta)}  \Biggr\}.
\end{eqnarray*}
Due to assumption (\ref{thm:consist:assumption}), $\sqrt n L_n^{21}(\cdot)$ is bounded by the following variable,
which is tight:
\[
2\frac{\sigma^*}{J} \sum_{1\leq|l|\leq m_n}|c_l(f^*)|\sum_{j=1}^J
 | \xi_{j,l} |.
\]
Thus, $L_n^{21}$ converges uniformly in probability to $0$.
Similarly, $L_n^{22}$ is bounded by
\[
L_n^{2B}=\frac{1}{\sqrt n} \Biggl(\sum_{|2p|>n} |c_p(f^*)| \Biggr)\sum
_{|l|\leq m_n}\sum_{j=1}^J  | \xi_{j,l} |.
\]
Consequently, from assumption (\ref{thm:consist:assumption}),
$L_n^{22}$ converges uniformly in probability to~$0$.
Therefore, $L_n$ converges uniformly in probability to $0$.

It remains to prove that $D_n$ converges uniformly to $M$.
First it is easy to prove that $D_n^1$ converges to $D^1$ and $D_n^2$
converges to $D^2$, where
\begin{eqnarray}
\nonumber D^1(\beta)&=&\frac{1}{J}\sum_{j=1}^J\int_0^{2\pi}
\bigl(f_j^*(t)-\nu_j \bigr)^2  \frac{dt}{2\pi}\quad\mbox{and}
\\
\nonumber D^2(\beta)
&=&\sum_{l\in\mathbb{Z}^*}  |c_l(f^*)\phi(l\theta - l\theta^* ,a) |^2 .
\end{eqnarray}
Consequently, $D_n$ pointwise converges to $M=D^1-D^2.$ We prove now
that the convergence is uniform.
For all $\beta\in\mathcal{A}_0$, we have
%
\begin{eqnarray}
\nonumber|D_n^1 - D^1 | (\beta)
&\leq&
\frac{1}{J} \sum_{j=1}^J
\Biggl\{\Biggl| \int_0^{2\pi}{f_j^*(t)}^2\frac{dt}{2\pi}
- \frac{1}{n} \sum_{i=1}^n {f^*(t_i)}^2 \Biggl|
\\
&&\hspace*{37pt}{}+ 2 \upsilon_{\max} \Biggr|c_0(f_j^*) - \frac{1}{n} \sum_{i=1}^n f^*(t_i) \Biggl|
\Biggl\},
\nonumber
\\
\nonumber
|D_n^2-D^2 | (\beta)
&\leq&
\sum_{|l|>m_n}|c_l(f^*)|^2+ | D_n^{2B}(\beta) |
\\
\eqntext{\mbox{where } D_n^{2B} = 2\displaystyle\sum_{1\leq|l|<m} \Re
\{c_l(f^*)\phi(l\theta-l\theta^*,a)\overline{g_n^l(\beta)} \}
+\displaystyle\sum_{1\leq|l|<m}  |g_n^l(\beta) |^2.}
\end{eqnarray}
Using the Cauchy--Schwarz inequality and inequality (\ref{proof:borne:gnl}), we have that
\begin{eqnarray*}
 |D_n^{2B}(\beta) |\leq2\sum_{|l|<m}|c_l(f^*)|\sum_{|p|>m_n}|c_p(f^*)| +2m_n \Biggl|\sum_{|p|>m_n}|c_p(f^*)| \Biggr|^2.
\end{eqnarray*}
The assumption (\ref{thm:consist:assumption}) ensures the uniform
convergence of $D_n^{2B}$. Consequently, since $f^*$ is continuous, we
deduce the uniform convergence of $D_n^1$ and $D_n^2.$
\end{pf}

\begin{lemma}[(Uniqueness of minimum)] \label{lem:minimum}
$M$ has a unique minimum reached in point $\beta=\beta^*$.
\end{lemma}

\begin{pf}
First, $M^2$, $M^1$ are nonnegative functions and we have that \mbox{$M(\beta^*)=0$.}
Consequently, the minimum of $M$ is reached in $\beta=(\theta
,a,\upsilon)\in\mathcal{A}_0$ if and only if $M^1(\beta)=M^2(\beta)=0.$

But if $M^2$ is equal to $0,$ this implies that $\upsilon=\upsilon^*.$

Furthermore, using the Cauchy--Schwarz inequality, we have for all
$l\in\mathbb{Z}^*$ that $|\phi(l\theta,a)| \leq  1.$ Since there exist
$l\in\mathbb{Z}^*$ such that $c_l(f^*)\neq0$ ($f^*$ is not
constant), $M^1$
is equal to $0$ if and only if the vectors $(a_j^*)_{j=1,\ldots, J}$ and
$(a_j e^{il(\theta_j-\theta_j^*)})_{j=1,\ldots, J}$ are proportional
for such $l$.
From the identifiability constraints on the model, we deduce that
\begin{eqnarray*}
a=a^*\mbox{ and }\forall l \in\mathbb{Z}\mbox{ such }|c_l(f)|
\neq  0\qquad l(\theta^* - \theta)\equiv0\ (2\pi).
\end{eqnarray*}
Thus it suffices that $c_1(f)\neq0,$ or there exist two relatively
prime integers $l,k$ such that $c_l(f^*)\neq0, c_k(f^*)\neq0$ in
order that $\theta=\theta^*.$
In other words, 2$\pi$ is the minimal period of the function $f^*.$ In
conclusion, $M^1(\beta)$ is equal to zero if and only if $a=a^*$ and
$\theta=\theta^*.$
\end{pf}

\subsection{\texorpdfstring{Proof of Proposition
\protect\ref{prop:LAN}}{Proof of Proposition 3.1}}
The proof is divided in two parts. First, we prove that
$ \langle\cdot,\cdot \rangle$ is an inner product. Next, we have to choose suitable
points $(\alpha_n(h),f_n(h))$ in order to establish the LAN property.

\subsubsection*{$ \langle\cdot,\cdot \rangle$ is an inner product in
$\mathcal H$} The form $\langle\cdot,\cdot\rangle_{\mathcal{H}}$ is bilinear,
symmetric and positive. In order to be an inner product, the form
$\langle\cdot,\cdot\rangle_{\mathcal{H}}$ has to be definite. In other words, if
$h\in\mathcal{H}$ is such that $\|h\|_{\mathcal{H}}=0,$ we want to
prove that $h=0.$ Let $h$ be such a vector; then we have that
$h_\sigma=0$ and for all $j=2,\ldots, J$,
\begin{eqnarray}\label{proof:LAN:definite}
\|a_j^*h_f+h_{a,j}f^*-h_{\theta,j}a_j^*\partial f^*+h_{\upsilon,j}
\|_{\mathbb{L}^2}&=&0\quad\mbox{and}
\nonumber\\[-8pt]\\[-8pt]
\nonumber\Biggl\|a_1^*h_f-\frac{\rho}{a_1^*}f^*+h_{\upsilon,1}
\Biggr\|_{\mathbb{L}^2}&=&0,
\end{eqnarray}
where $\rho=\sum_{k=2}^Jh_{a,k} a_k^*.$
Since the functions $h_f$, $f^*$ and $\partial f^*$ are orthogonal to
$1$ in $\mathbb L^2[0,2\pi]$, we deduce that $h_{\upsilon,j}=0$ for
all $j.$
Moreover, the functions $h_f$ and $f^*$ are continuous and the equation
(\ref{proof:LAN:definite}) implies that ${a_1^*}^2h_f=\rho f^*$ and
that for all $j=2,\ldots, J$ ($f^*$ and $\partial f^*$ are orthogonal),
\[
\biggl\| \biggl(\frac{a_j^*\rho}{{a_1^*}^2}+h_{a,j} \biggr)f^* \biggr\|
_{\mathbb{L}^2}   = 0\quad\mbox{and}\quad\|h_{\theta,j}a_j^*\partial f^*\|_{\mathbb{L}^2}  = 0.
\]
Since $f^*$ is not constant, we deduce that for all $j=2,\ldots, J$ that
$h_{\theta,j}=0$ and $a_j^*\rho/{{a_1^*}^2}+h_{a,j}=0.$ Consequently,
$\rho$ verifies the equation $\rho\frac{J-{a_1^*}^2}{{a_1^*}^2}+\rho
=0.$ Then $\rho$ is equal to zero and $h=0.$

\subsection*{The LAN property}
Let $h$ be in $\mathcal H$.
In order to satisfy the identifiability constraints of the model, we
choose the sequences $(\alpha_n(h),f_n(h))$ [with $\alpha_n(h)
= ((\theta_n^{(j)}(h))_{1\leq j\leq J} ,(a_n^{(j)}(h))_{1\leq j\leq J},
(\upsilon_n^{(j)}(h))_{1\leq j\leq J} ,\sigma_n(h))$]
such that
\begin{eqnarray*}
\theta_n^{(j)}(h)
&=& \theta_j^* + \frac{1}{\sqrt n}h_{\theta,j}\quad\mbox{and}\quad a_n^{(j)}(h) = a_j^* +
\frac{1}{\sqrt n}h_{a,j}\qquad\forall j=2,\ldots, J,
\\
\theta_n^{(1)}(h)
&=& 0,\qquad a_n^{(1)}(h) =  \sqrt{J - \sum_{j=2}^Ja_n^{(j)}(h)^2}\quad\mbox{and}\quad\sigma_n(h) =
\sigma^* + \frac{h_\sigma}{\sqrt n},
\\
f_n(h) &=& f_n = f^* + \frac{1}{\sqrt n}h_f\quad\mbox{and}\quad
\upsilon_n^{(j)}(h) = \upsilon_j^* + \frac{1}{\sqrt n}h_{\upsilon,j}\qquad\forall j=1,\ldots, J.
\end{eqnarray*}
Using the uniform continuity of $\partial f^*$ and $h_f$, we uniformly
establish for $i=1,\ldots, n$ that
\begin{eqnarray*}
f_n\bigl(t_i-\theta_n^{(j)}(h)\bigr)-f_n(t_i-\theta_j^*)
&=&\frac{h_{\theta,j}}{\sqrt n}\partial f^*(t_i-\theta_j^*)+o\bigl(1/\sqrt n\bigr)\quad\forall j=1,\ldots,
J,
\\
\bigl(a_n^{(1)}(h)-a_1^*\bigr) f_n^* (t_i )
&=&-\frac{\sum_{j=2}^Jh_{a,j} a_j^*}{a_1^*\sqrt n}f^* (t_i
)+o\bigl(1/\sqrt{n}\bigr),
\\
\log \biggl( 1+\frac{h_\sigma/\sigma^*}{\sqrt n} \biggr)
&=&\frac{h_\sigma/\sigma^*}{\sqrt n}-\frac{(h_\sigma/\sigma^*)^2}{n}+o(n^{-1}).
\end{eqnarray*}
Then, with the notation of the proposition, we may deduce that
\[
\Lambda_n(\alpha_n(h),f_n(h)) = \Delta_{n}(h) - \frac
{1}{2n}\sum_{i=1}^n \sum_{j=1}^J  {A_{i,j}^n(h)}^2  - \frac
{J}{2}\frac{\sigma^2}{{\sigma^*}^2} + o_{\mathbb{P}}(1).
\]
$\sum_{i=1}^n \sum_{j=1}^J  {A_{i,j}^n( h )}^2 / n $ is a
Riemann sum which converges to $\|h\|^2_{\mathcal{H}}$.
Moreover, from the Lindeberg--Feller central limit theorem (see
\cite{VandVaart}, Chapter~2)
$\Delta_{n}(h)$ converges in distribution
to $\mathcal{N} (0,\|h\|^2_{\mathcal{H}} )$.

\subsection[The efficient estimation]{The efficient estimation of
$\theta^*,$ $a^*$ and $\upsilon^*$}
\begin{lemma}[(The derivative of $\nu$)]\label{lem:derivative:nu}
The representant of the $\nu_n$'s derivative is
$\dot{\nu}=((\dot{\nu}^{\theta}_j)_{2\leq j\leq J},(\dot{\nu}^{a}_j)_{2\leq j\leq J},
(\dot{\nu}^{\upsilon}_j)_{1\leq j\leq J}\in\overline{\mathcal{H}}^{3J-2},$ where
\begin{eqnarray*}
\dot{\nu}^{\theta}_j
&=&\frac{{\sigma^*}^2}{\|\partial f^*\|_{\mathbb{L}^2}}
\biggl(\dot{\theta}^j,0,0,0,\frac{1}{{a_1^*}^2}\partial f^* \biggr)\qquad\mbox{for}  j= 2,\ldots,
J,
\\
\dot{\nu}^{a}_j
&=&\frac{{\sigma^*}^2}{\|f^*\|^2_{\mathbb{L}^2}}
(0,\dot{a}^j,0,0,0 )\hspace*{62pt}\mbox{for}  j=2,\ldots, J,
\\
\dot{\nu}^{\upsilon}_j
&=& (0,0,e_j,0,0 )\hspace*{98pt}\mbox{for}  j= 1,\ldots, J,
\end{eqnarray*}
where the vector $e_j$ is the $j$th vector of canonical basis of
$\mathbb{R}^J$, and the vectors $\dot{\theta}^j=(\dot{\theta}^j_k)_{k=2,\ldots, J}$
and $\dot{a}^j=(\dot{a}^j_k)_{k=2,\ldots, J}$ are defined as
\begin{eqnarray*}
\dot{\theta}_{k}^j&=&
\cases{
1/{a_1^*}^2, &\quad if $k\neq j$,\cr
1/{a_1^*}^2+1/{a_j^*}^2, &\quad if $k=j$,
}
\\
\dot{a}_{k}^j&=&
\cases{
-a_2^*a_k^*/J,&\quad if $k\neq j$,\cr
1-{a_k^*}^2/J,&\quad if $k=j$.
}
\end{eqnarray*}
\end{lemma}

\begin{pf}
For $h\in\mathcal H$ and $h'\in\mathcal H$, we may rewrite the inner
product of the tangent space under the following form:
\begin{eqnarray*}
{\sigma^*}^2 \langle h,h'\rangle
&=& Jh_\sigma h_\sigma'+\langle h_f,Jh_f'
-\lambda\partial f^*\rangle
\\
&&{}+\sum_{k=2}^Jh_{\theta,k}
\langle\partial f^*,-a_k^*h_f'+h_{\theta,k}
a_k^*\partial f^* \rangle
\\
&&{}+\sum_{k=2}^J h_{a,k} \biggl\langle f^*,h_{a,k}' f^*+\frac{a_k^*}{{a_1^*}^2}\rho f^*\biggr\rangle
+\sum_{k=2}^Jh_{\upsilon,k}h_{\upsilon,k}',
\end{eqnarray*}
where $\lambda=\sum_{k=2}^Jh_{\theta,k}' a_k^*$ and $\rho=\sum_{k=2}^J h_{a,k}' a_k^*$.
Let $k\in\{2,\ldots, J\}$ be a fixed integer;
we want to find $h'$ such that for all $h\in\mathcal{H},$
$\langle h,h'\rangle=h_{\theta,k}.$ Consequently, such~$h'$ verifies
these equations:
%
\begin{eqnarray}
h_{f}=\lambda\partial f^*/J,\qquad h_{\sigma}'=0\quad\mbox{and}\quad
h_{\upsilon,j}'
&=& 0,\qquad\forall j=1,\ldots, J,
\\
(h_{a,j}'+\rho a_j^*/{a_1^*}^2 )\|f^*\|^2
&=& 0,\qquad\forall
j=2,\ldots, J,\label{proof:derivative:nu:a}
\\
(-\lambda/J+h_{\theta,j} )\|\partial f^*\|^2
&=&
\cases{
{\sigma^*}^2,&\quad if $j=k$,\cr
0,&\quad if $j\neq k$.
}\label{proof:derivative:nu:theta}
\end{eqnarray}
Combining equations (\ref{proof:derivative:nu:a}) and (\ref{proof:derivative:nu:theta}), we have that
\begin{eqnarray*}
\lambda{a_1^*}^2\|\partial f^*\|^2/J={\sigma^*}^2\quad\mbox{and}\quad\rho J\| f^*\|^2/{a_1^*}^2=0.
\end{eqnarray*}
Thus we deduce that $\rho=0$ and $\lambda=J{\sigma^*}^2/
({a_1^*}^2\|\partial f^*\|^2 )$.
Consequently, $h'$ is equal to $\dot{\nu}^{\theta}_k$.

We likewise solve the equation $ \langle h,h'\rangle=h_{a,k}.$
Finally, we have that $\|f^*\|^2 \rho={\sigma^*}^2 {a_1^*}^2 a_k^*/J$
and $\lambda=0.$ Hence the solution is $h'=\dot{\nu}^{a}_k$.
\end{pf}

\begin{prop}\label{prop:gradient}
Under the assumptions and notation of Theorem \ref{thm:efficiency}, we
have that
\[
\sqrt{n} \nabla M_n(\beta^*)=-\frac{2\sigma^*}{J}G_n+o_{\mathbb P}(1)\qquad\mbox{where}\quad{}^t\!G_n
={}^t\!(G_n^\theta, G_n^a, G_n^\upsilon ).
\]
$G_n$ is a Gaussian vector which converges in distribution to $\mathcal{N}_{3J-2} (0,H)$ and is defined as
\begin{eqnarray*}
G_n^\theta
&=&\frac{1}{\sqrt n}\sum_{i=1}^n
\biggl[
\frac{a_1^*}{J}A^2 \mano{-}D+\frac{1}{J}A^2\,\,{}^t\!\!A
\biggr]\partial F^*(t_i)\varepsilon_{i,\bolds{\cdot}},
\\
G_n^a
&=&
\frac{1}{\sqrt n}\sum_{i=1}^n
\biggl[
\frac{-1}{a_1^*}A \mano I_{J-1}
\biggr]F^*(t_i)\varepsilon_{i,\bolds{\cdot}},
\\
G_n^\upsilon
&=&\frac{1}{\sqrt n}\sum_{i=1}^n
\varepsilon_{i,\bolds{\cdot}}\quad\mbox{and}\quad{}^t
\varepsilon_{i,\bolds{\cdot}}={}^t (\varepsilon_{i,1},\ldots,
\varepsilon_{i,J})\qquad\mbox{for } i=1,\ldots, n.
\end{eqnarray*}
\end{prop}

\begin{pf}
In order to prove that proposition, we proceed in two steps. First,
using the notation of Proposition \ref{lem:contraste}, we show that
\[
\sqrt n \nabla M_n(\beta^*)
=\sqrt n \bigl(\nabla L^1_n(\beta^*)
-L_n^{21}(\beta^*) \bigr)=-\frac{2\sigma^*}{J}G_n+o_{\mathbb P}(1).
\]
At the end, we prove that $(G_n^\theta ,G_n^a ,G_n^\upsilon )$ is
a Gaussian vector which converges to $\mathcal{N}_{3J - 2}(0,H )$.

First, we study singly the gradient of $G_n$, $L_n$ and $Q_n$.
Let $k\in \{2,\ldots, J \}$ be fixed. The partial
derivative with respect to the variable $\theta_k$ is
\begin{eqnarray*}
\frac{\partial Q_n}{\partial\theta_k}(\beta^*)
&=&-\frac{2}{n}
\sum_{1\leq|l|<m_n}\Re \biggl( \frac{ila_k^*e^{il\theta_k^*}}{J}\xi_{k,l}\overline{\xi_l(\beta^*)}  \biggr).
\end{eqnarray*}
It is bounded by
\begin{eqnarray*}
\Biggl|\sqrt n \frac{\partial Q_n}{\partial\theta_k}(\beta^*)\Biggr|
\leq\frac{2}{J^2\sqrt n}\sum_{1\leq|l|<m_n} |l|  |\xi_{k,l} |\sum_{j=1}^J |\xi_{j,l} |.
\end{eqnarray*}
Thus $\sqrt n \frac{\partial Q_n}{\partial\theta_j}(\beta^*)$
converges in probability to $0$ if $m_n^4/n=o(1).$ Similarly, the
partial derivative with respect to the variable $a_k$ converges in
probability to $0$, too. Consequently, $\sqrt n \nabla Q_n(\beta^*)$
converges to $0$ in probability.

Concerning the deterministic part, the partial derivative with respect
to $\theta_k$ is\looseness=1
\begin{eqnarray*}
&&\frac{\partial D_n}{\partial\theta_k} ( \beta^* )
=  -\frac{2}{J}\sum_{1\leq|l|<m_n}\Re
\Biggl\{il\sum_{p\in\mathbb{Z}}  c_p(f^*)\varphi_n\biggl(\frac{l-p}{n}  \biggr)
{a_k^*}^2e^{i(l-p)\theta_k^*}
\\
&&\hspace*{125pt}{}\times\Biggl(\sum_{p\in\mathbb{Z}}  \overline{c_p(f^*)}\varphi_n
\biggl(\frac{p-l}{n}\biggr)\varphi\bigl((p - l)\theta^* , a^*\bigr)\Biggr)
\Biggr\}.
\end{eqnarray*}
Using the inequality (\ref{proof:borne:gnl}), it is bounded by
\begin{eqnarray*}
&&\sqrt n \frac{\partial D_n}{\partial\theta_k}(\beta^*)=\frac{2{a_k^*}^2}{J}
\sqrt n  \Biggl\{ 2 \sum_{|l|\leq m_n}|lc_l(f^*)| \sum_{2|p|\geq n}
|c_p(f^*)|
\\
&&\qquad{}\hspace*{102pt}+ \sum_{|l|\leq m_n}|l|
 \Biggl(\sum_{2|p|\geq n} |c_p(f^*)|\Biggr)^2 \Biggr\}.
\end{eqnarray*}
Consequently, we deduce from the assumptions of the theorem that $\sqrt
n \frac{\partial D_n}{\partial\theta_k}(\beta^*)$ converges in
probability to $0.$
In like manner, $\sqrt n \frac{\partial D_n}{\partial a_k}(\beta^*)$
converges in probability to $0$, too.
For the partial derivative with respect to $\upsilon_k$, we have
\begin{eqnarray*}
\sqrt n\frac{\partial D_n}{\partial\upsilon_k}(\beta^*) = -\frac
{2\sqrt n}{Jn} \Biggl(\sum_{i=1}^nf^*_k(t_i) - \upsilon_k^* \Biggr)
=  -\frac{2a_k^*}{J} \sqrt n \sum_{p\in n\mathbb
{Z}^*}c_p(f^*)e^{-ip\theta_k^*}.
\end{eqnarray*}
Thus from assumption (\ref{thm:efficiency:asump:m}), we deduce that
$\sqrt n\partial D_n/\partial\nu_k$ in $\beta^*$ converges to~$0$.
Finally, $\sqrt n \nabla D_n(\beta^*)$ converges to $0$ in probability.

Therefore, we have that $\sqrt n \nabla M_n(\beta^*)=\sqrt n \nabla
L_n(\beta^*)+o_{\mathbb{P}}(1).$
With the notation of Lemma \ref{lem:contraste}, we have
\begin{eqnarray*}
&&\sqrt n \frac{\partial L_n^{22}}{\partial\theta_k}(\beta^*) =
\frac{2}{J}\sum_{1\leq|l|<m_n}  \Re \Biggl\{il a_k^*
\bigl(-e^{il\theta_k^*}\xi_{k,l}+ a_k^*\bar{\xi_l}(\beta) \bigr)
\mathop{\sum_{2|p|\geq n}}_{p-l \in n\mathbb{Z}}c_p(f^*) \Biggr\}.
\end{eqnarray*}
The centered Gaussian variable $\sqrt n \frac{L_n^{22}}{d \theta_k}(\beta^*)$ has a variance bounded by
\[
 \Biggl( \sum_{2|p|>n} |c_p(f^*)| \Biggr)^2   2m_n^3.
\]
From assumption (\ref{thm:efficiency:asump:m}), we conclude that
$\sqrt n \frac{\partial L_n^{22}}{d \theta_k}(\beta^*)$ converges to
$0$ in probability.
In like manner, $\sqrt n \frac{\partial L_n^{22}}{\partial a_k}(\beta^*)$ converges in probability to $0$, too.
Thus we have that
$\sqrt n \nabla M_n(\beta^*)=\sqrt n \nabla L_n^1(\beta^*)-\sqrt n \nabla L_n^{21}(\beta^*)+o_{\mathbb{P}}(1).$
After straightforward computations, we obtain
\begin{eqnarray*}
\sqrt n\frac{\partial M_n}{\partial\theta_k}(\beta^*)
&=&-\frac{2\sigma^*}{J}\sum_{1\leq|l|\leq m_n}\Re \bigl\{lc_l(f^*)
\bigl(a_k^*\overline{\xi_l(\beta^*)}-a_k^*e^{-il\theta_k^*}\overline{\xi_{k,l}}
\bigr)\bigr\},
\\
\sqrt n\frac{\partial M_n}{\partial a_k}(\beta^*)
&=&-\frac{2\sigma^*}{J}\sum_{1\leq|l|\leq m_n}\Re \biggl\{c_l(f^*) \biggl(e^{-il\theta_k^*}
\overline{\xi_{k,l}}-\frac{a_k^*}{a_1^*}\overline{\xi_{1,l}} \biggr)
\biggr\},
\\
\sqrt n\frac{\partial M_n}{\partial\upsilon_k}(\beta^*)
&=&-\frac{2\sigma^*}{J\sqrt n}\sum_{i=1}^n\varepsilon_{i,k}=-\frac{2\sigma^*}{J}\overline{\xi_{k,0}}.
\end{eqnarray*}
We can now define $(G_n^\theta,G_n^a,G_n^\upsilon)$ as
\begin{eqnarray*}
G_n^\theta
&=&
\sum_{1\leq|l|\leq m_n} \Re \biggl\{ilc_l(f^*)\biggl(\frac{a_1^*}{J}
\mathbb A^2 \mano{-}D+\frac{1}{J}\mathbb{A}^2\,\,{}^t\!\mathbb{A} \biggr)X_l^*
\biggr\}+o_{\mathbb{P}}(1),
\\
G_n^a
&=& \sum_{1\leq|l|\leq m_n} \Re \biggl\{c_l(f^*)
\biggl(\frac{-1}{a_1^*}\mathbb A \mano I_{J-1}\biggr)
X_l^* \biggr\}+o_{\mathbb{P}}(1)\quad\mbox{and}
\\
G_n^\upsilon
&=& \Re \{X_0^* \}+o_{\mathbb{P}}(1),
\end{eqnarray*}
where $X_l^*$ denote the independent identically distributed complex
Gaussian vectors defined as
\[
{}^t\! X_l^*={}^t\!(e^{-il\theta_k^*}\overline{\xi_{1,l}},
\ldots, e^{-il\theta_k^*}\overline{\xi_{J,l}} ) .
\]

Since $G_n^\theta$ and $G_n^a$ do not depend on $X_0^*,$ $G_n^\upsilon
$ is independent of $G_n^\theta$ and $G_n^a.$ Moreover, its variance
matrix is equal to the identity matrix of $\mathbb{R}^J.$
Furthermore, the imaginary part and the real part of $c_l(f^*)X_l^ *$
are independent. Consequently, $G_n^\theta$ and $G_n^a$ are
asymptotically independent with covariance matrix
$\|\partial f^*\|^2 \times(D^2-A^2\,\,{}^t\!\!A^2/J )$ and
$\| f^*\|^2 (I_{J-1}- A\,\,{}^t\!\!A/{a_1^*}^2 )$, respectively.

By the definition of $(\xi_{k,l})$ (Remark \ref{remqnotations}), we
deduce from assumption (\ref{thm:efficiency:asump:f}) that for a
fixed $k=1,\ldots, J,$
\begin{eqnarray*}
\Re \Biggl\{\sum_{|l|\leq m_n} il c_l(f^*)\bar{\xi}_{k,l}e^{-il\theta_k^*} \Biggr\}
&=&\frac{1}{\sqrt n}\sum_{i=1}^n
\varepsilon_{i,k}\Re \Biggl\{\sum_{|l|\leq m_n} il c_l(f^*)e^{il(t_i-\theta_k^*)} \Biggr\}
\\
&=&\frac{1}{\sqrt n}\sum_{i=1}^n \varepsilon_{i,k}\partial f^*(t_i-\theta_k^*)+o_{\mathbb{P}}(1) .
\end{eqnarray*}
Thus, $(G_n^\theta,G_n^a,G_n^\upsilon)$ are equal to the expression
defined in the proposition.\looseness=1
\end{pf}

\begin{prop}\label{prop:hessian}
Under the assumptions and notation of Theorem \ref{thm:efficiency}, we have
\[
\nabla^2 M_n(\bar{\beta}_n)
\mathop{\stackrel{\mathbb{P}_{\beta^*}}{\longrightarrow}}_{n\rightarrow\infty}
 -\frac{2}{J^2}H.
\]
\end{prop}

\begin{pf}
The matrix $-2H/J^2$ is the value of the Hessian matrix of $M$ in
point~$\beta^*$.
We study locally the Hessian matrix of $M_n$. Consequently, we may
assume that the sequences ($\bar{\beta}_n$) are in the following set:
\[
A_n^{\mathit{loc}}= \{ (\theta,a,\upsilon)\in\mathcal A_0,  a_1>r
\mbox{ and }\|\beta -  \beta^* \|\leq\|\hat{\beta}_n  -
\beta^* \| \},
\]
where $a_1^*>r>0.$ Notice that for $\varepsilon>0,$ we have
\begin{eqnarray*}
&&\mathbb{P} \biggl(\sup_{\beta\in A_n^{loc}}  \|\nabla^2
M_n(\beta) - \nabla^2 M(\beta^*) \|>2\varepsilon \biggr)
\\
&&\qquad\leq\mathbb{P} \biggl(\sup_{\beta\in A_n^{loc}}  \|\nabla^2
M_n(\beta) - \nabla^2 M(\beta) \|>\varepsilon \biggr)
\\
&&\qquad\quad{}+ \mathbb{P} \biggl(\sup_{\beta\in A_n^{loc}}  \|\nabla^2
M(\beta)-\nabla^2 M(\beta^*) \|>\varepsilon \biggr).
\end{eqnarray*}
As in Lemma \ref{lem:contraste}, assumptions (\ref
{thm:efficiency:asump:f}) and (\ref{thm:efficiency:asump:m}) assure
the uniform convergence in probability of $\nabla^2 M_n$ to the
Hessian matrix of $M$ on $A_n^{loc}$.
Thus, the first term of inequality converges to $0$ with $n$.

Since $\nabla^2 M$ is continuous in $\beta^*$, there exists $\delta
>0$ such that
\[
\nabla^2 M (B(\beta^*,\delta) )\subseteq B (\nabla
^2 M(\beta^*),\varepsilon ).
\]
Consequently, we have the following inclusion of event:
\[
\biggl (\sup_{\beta\in A_n^{loc}}  \|\nabla^2 M(\beta)-\nabla^2 M(\beta^*) \|>\varepsilon \biggr) \subseteq (  \|
\hat{\beta}_n  -  \beta^*  \|> \delta ).
\]
Thus, from Theorem \ref{thm:consistance}, the second term of the
inequality converges to $0$, too.
\end{pf}

\subsubsection{The estimation of the common shape}
\begin{remq}\label{remq:smooth:f}
If the assumptions of Theorem \ref{thm:efficiency} hold, we obtain
using the Cauchy--Schwarz inequality that
\[
\sum_{|l|>n}|c_l(f^*)|\leq \Biggl\{\sum_{|l|>n}|lc_l(f^*)| \Biggr\}
^{1/2} \Biggl\{\sum_{|l|>n}|lc_l(f^*)|/l^2 \Biggr\}^{1/2}=o(1/n).
\]
Similarly, if $f^*$ is $k$ times differentiable and $f^{(k)}$ is
squared integrable, we have
\[
\sum_{|l|>n}|c_l(f^*)|=o(n^{-k+1/2})\quad\mbox{and}\quad\sum
_{|l|>n}|c_l(f^*)|^2=o(n^{-2k}).
\]
\end{remq}
\begin{pf*}{Proof of Corollary \ref{cor:shape:rate}}
Using the notation of Lemma \ref{lem:contraste}, we have for all
$t\in\mathbb R,$
%
\begin{eqnarray}
\hspace*{25pt}&&{f}^*(t) - f_{{\beta}^*}(t)
\nonumber\\
&&\qquad= \sum_{|l|> m_n}c_l(f^*)e^{ilt}+\sum_{1\leq|l|\leq m_n}e^{ilt}
\sum_{|2p|> n, p-l \in n\mathbb{Z}}  c_p(f^*)\phi(l\hat\theta-p\theta^*,\hat a)\nonumber
\\
&&\qquad\quad{}+ \sum_{1\leq|l|\leq m_n} c_l(f^*)\bigl\{\phi\bigl(l(\hat\theta-\theta^*),
\hat a\bigr)-1 \bigr\}e^{ilt}\label{function:term2}
\\
&&\qquad\quad{}+\sigma^*\sum_{1\leq|l|\leq m_n}\xi_l(\hat\beta)\frac{e^{ilt}}{\sqrt n}\label{function:term3}.
\end{eqnarray}
Since Theorem \ref{thm:efficiency} holds and using the delta method,
we have for all $j=1,\ldots, J,$
\[
e^{il(\hat{\theta}_j-\theta_j^*)}-1=il(\hat{\theta}_j-\theta
_j^*)+o_{\mathbb P}\bigl(l/\sqrt n\bigr).
\]
Moreover, we have
\[
 \bigl|\phi\bigl(l(\hat{\theta}-\theta^*),\hat{a}\bigr)-1 \bigr|\leq\frac
{1}{J}\sum_{j=1}^Ja_j^*|\hat{a}_{j}-a_j^*|+\frac{1}{J}\sum
_{j=1}^J{a_j^*}^2\bigl|e^{il(\hat{\theta}_j-\theta^*_j)}-1\bigr|.
\]
Then, we deduce that
\begin{eqnarray*}
\sup_{t\in\mathbb{R}}|(\ref{function:term2})|
=\mathcal{O}_{\mathbb P}\bigl(1/\sqrt n\bigr)\quad\mbox{and}\quad\mathbb E
\|(\ref{function:term2})\|_{\mathbb{L}^2}^2=\mathcal{O}(1/n).
\end{eqnarray*}
Using (\ref{proof:borne:xil}) and the Cauchy--Schwarz inequality, we have
\begin{eqnarray*}
W_n&=&\Biggl|\sum_{1\leq|l|\leq m_n}\xi_l(\hat\beta)
\frac{e^{ilt}}{\sqrt n} \Biggr| \leq\frac{1}{\sqrt{Jn}}\sum_{1\leq|l|\leq
m_n}\sum_{j=1}^J|\xi_{j,l}|,
\\
\int_0^{2\pi}\mathbb E W_n^2\frac{dt}{2\pi}
&=&\frac{1}{nJ}\sum_{1\leq|l|\leq m_n}\sum_{j=1}^J|\xi_{j,l}|^2.
\end{eqnarray*}
Hence we deduce by the Markov inequality that
\[
W_n=O_{\mathbb P} \bigl(m_n/\sqrt n \bigr)\quad\mbox{and}\quad\int_0^{2\pi}
\mathbb E W_n^2\frac{dt}{2\pi}=O(m_n/n).
\]
Then, using Remark \ref{remqnotations}, the corollary results.
\end{pf*}

\begin{lemma}\label{lem:asymptotic:linear:f}
Let $l$ be in $\mathbb{Z}^*.$
For a large $n$, we have
\begin{eqnarray*}
\sqrt n \Re \bigl(\hat{c}_l(\hat\beta_n)-c_l(f^*) \bigr)
&=&\Delta_n (-\Re(lc_l(f^*)\tilde{h}^f )+\Delta_n
\biggl(0,0,0,0,\frac{\cos(l\cdot)}{J} \biggr)
\\
&&{}+o_{\mathbb P}(1),
\\
\sqrt n \Im \bigl(\hat{c}_l(\hat\beta_n)-c_l(f^*) \bigr)
&=&\Delta_n (\Im(lc_l(f^*)\tilde{h}^f )+\Delta_n\biggl(0,0,0,0,\frac{-\sin(l\cdot)}{J}
\biggr)
\\
&&{}+o_{\mathbb P}(1),
\end{eqnarray*}
where $\tilde{h}^f=\frac{{\sigma^*}^2}{\| \partial f^*\|_{\mathbb{L}^2}} (\frac{\mathbb{I}_{J-1}}{{a_1^*}^2},0,0,0,\frac
{J-{a_1^*}^2}{{a_1^*}^2}\partial f^* ).$
\end{lemma}

\begin{pf}
Let $l$ be in $\mathbb{Z}^*.$
For $n$ large enough (such wise $|l|\leq m_n$), from the continuous
mapping theorem \cite{VandVaart}, Theorem 2.3, and from assumption (\ref{thm:efficiency:asump:f})
ensures that $c_l^n(f^*)$ converges to
$c_l(f^*)$ with a speed $\sqrt n,$ we obtain
\begin{eqnarray*}
\sqrt n c_l  (\hat{f}_n-f^* )
&=&\sqrt n  \bigl(\hat
{c}_l(\hat\beta_n)-c_l(f^*) \bigr),
\\
&=&c_l(f^*)\sqrt n  \Biggl(\frac{1}{J}\sum_{j=1}^J\hat
{a}_{j,n}a_j^*e^{il(\hat\theta_{j,n}-\theta_j^*)}-1 \Biggr)
+\xi_l(\theta^*,a^*)+o_{\mathbb P}(1).
\end{eqnarray*}
Since $\sqrt n(\hat\theta_n - \theta^* ,\hat{a}_n - a^*)$
converges in distribution (Theorem \ref{thm:efficiency}),
we use the delta method (\cite{VandVaart}, Chapter~3):
\begin{eqnarray*}
\sqrt n c_l  (\hat{f}_n-f^* )
=ilc_l(f^*)\sum_{j=2}^J\frac{{a_j^*}^2}{J}\sqrt n(\hat\theta_j-\theta_j^*)+\xi_l(\theta^*,a^*)+o_{\mathbb P}(1).
\end{eqnarray*}
Thus from Theorem \ref{thm:efficiency} and Lemma \ref
{lem:derivative:nu} and due to the linearity of $\Delta_n(\cdot),$ we have
\begin{eqnarray*}
\sqrt n c_l  (\hat{f}_n-f^* )
&=&ilc_l(f^*)\sum_{j=2}^J\frac{{a_j^*}^2}{J}\Delta_n(\tilde{h}^{\theta}_j)
+\xi_l(\theta^*,a^*)+o_{\mathbb P}(1)
\\
&=&\frac{ilc_l(f^*){\sigma^*}^2}{\| \partial f^*\|_{\mathbb{L}^2}}
\Delta_n (\tilde{h}^f )+\xi_l(\theta^*,a^*)+o_{\mathbb P}(1).
\end{eqnarray*}
Using the definition of $\xi_l$ (see Remark \ref{remqnotations}), we have
\begin{eqnarray*}
\hspace*{-12pt}\Re  ( {\xi_l(\theta^*,a^*)} ) =  \Delta_n  \biggl(0, \frac{\cos(l\cdot)}{J}
\biggr)\quad\mbox{and}\quad\Im({\xi_l(\theta^*,a^*)} ) =  \Delta_n \biggl(0, \frac{-\sin(l\cdot)}{J}  \biggr).
\end{eqnarray*}
\upqed
\end{pf}

\section*{Acknowledgments}\label{sec:acknow}
I am grateful to Professor Fabrice Gamboa for his advice during the preparation
of this article. I also would like to thank an Editor and the referees
for their
comments and suggestions that helped to greatly improve the
presentation of the paper.
And I gratefully acknowledge the webmaster of \url{http://www.engr.udayton.edu/weather}
for the permission to use the weather database.

\printaddresses

\end{document}